%% file: symcubic.tex
\documentclass[11pt]{amsart}
\usepackage{graphicx}
\usepackage{amssymb}
\usepackage{array}
\usepackage{amsthm}
\usepackage{mathtools}
\usepackage{microtype}
\usepackage{mathrsfs}
\usepackage[all]{xy}
\usepackage{fullpage}  
\usepackage{setspace}
\usepackage{comment}
\usepackage{tikz-cd}
\usepackage{dynkin-diagrams}
\usepackage{amsmath}
\usetikzlibrary{positioning}
\usepackage{amsthm,thmtools}
\usepackage[colorlinks,citecolor=darkgreen,linkcolor=darkgreen,urlcolor=darkblue]{hyperref}

\usepackage{float}

\input{commands.tex}
\input{thms.tex}

\providecommand{\C}{\mathbb{C}}
\renewcommand{\P}{\mathbb{P}}
\providecommand{\sm}{\mathrm{sm}}
\newcommand{\quot}{\backslash\!\!\backslash}

\newcommand{\sqrtone}{\alpha}
\newcommand{\sqrttwo}{\beta}

\title{Monodromy in the space of symmetric cubic surfaces with a line}
\author{Thomas Brazelton and Sidhanth Raman}
\begin{document}
\begin{abstract}
    We explore the enumerative problem of finding lines on cubic surfaces defined by symmetric polynomials. We prove that the moduli space of symmetric cubic surfaces is an arithmetic quotient of the complex hyperbolic line, and determine constraints on the monodromy group of lines on symmetric cubic surfaces arising from Hodge theory and geometry of the associated cover. This interestingly fails to pin down the entire Galois group. Leveraging computations in equivariant line geometry and homotopy continuation, we prove that the Galois group is the Klein 4-group. This means that, despite a general cubic surface admitting no formula in radicals for its lines, an $S_4$-symmetric cubic does; we work out these formulas explicitly. This is the first computation in what promises to be an interesting direction of research: studying monodromy in classical enumerative problems restricted by a finite group of symmetries.
\end{abstract}

\maketitle

\section{Introduction}

The \textit{Galois group} of an enumerative problem is a classical object of study in enumerative algebraic geometry. It was first introduced by Jordan as one of the main subjects of interest at the genesis of Galois theory \cite{Jordan}. This idea enjoyed a revival a century later when Harris proved that the Galois group of an enumerative problem agrees with the monodromy group of its associated cover \cite{Harris}. In modern mathematics, Galois groups can be approached from a wide number of perspectives, from Hodge theory and hyperbolic geometry \cite{AllcockCarlsonToledo}, to Lie theory \cite{Manivel}, to numerical analysis and homotopy continuation \cite{LeykinSottile}, to name a few.\footnote{For a lovely introduction to the history and appearance of Galois groups in enumerative geometry, we refer the reader to \cite{SottileYahl}.}

Contemporary geometers such as Klein were interested in exploring how symmetries of objects manifest in enumerating various quantities attached to them. Recent work of the first-named author introduces tools from equivariant homotopy theory to explore how Poncelet's principle of conservation of number interacts with symmetry, an example being that a smooth cubic surface defined by a symmetric polynomial always has the same $S_4$-symmetries on its lines \cite{EEG}. Such cubic surfaces are called \textit{symmetric cubic surfaces}.

In this paper we initiate an exploration of monodromy groups of symmetric enumerative problems. This flavor of question is well-studied in geometric group theory; for example, many have studied rigidity phenomena for finite index subgroups of lattices in Lie groups (e.g. \cite{Margulis} and \cite{FarbWeinberger}) and equivariant problems for their non-linear analogues like mapping class groups and $\mathrm{Out}(F_n)$ (e.g. \cite{BirmanHilden}, \cite{MaclachlanHarvey}, \cite{FarbHandel}, and \cite{LandesmanLittSawin}). However the setting we pursue is of a completely different shape --- since the Galois group of lines on a cubic surface (and many related problems) is finite, we cannot leverage such tools, e.g. Teichm\"uller theory, to approach this question, and alternative techniques are needed.

Our main result is a \textit{computation of the monodromy group of lines on symmetric cubic surfaces}, which we show is equal to the Klein $4$-group. This is carried out via a combination of moduli-theoretic techniques, classical analysis of the Weyl group of the $E_6$ lattice, as well as group-theoretic computations in GAP and contemporary certified tracking homotopy continuation algorithms. Along the way we prove that \textit{the moduli of stable symmetric cubic surfaces is an arithmetic quotient of the complex hyperbolic line}. The latter result mirrors the landmark work of Allcock, Carlson, and Toledo at the turn of the century \cite{AllcockCarlsonToledo}, where they show the moduli space of stable cubic surfaces is an arithmetic quotient of complex hyperbolic $4$-space. We explore the appearance of our Klein 4-group $K_4$ in both the Weyl group of $E_6$ and in the projective orthogonal group $\PO(4,1,\mathbb{F}_3)$. Finally we establish an \emph{explicit formula in radicals} for the 27 lines on an $S_4$-symmetric cubic surface. This is an interesting result because, famously, no formula in radicals exists for lines on a general cubic surface. The presence of symmetry drives down the size of the Galois group of such a problem, and converts an unsolvable problem into a solvable one.

\subsection{Main results}
Before we state our main theorems more formally, we fix some notation. Let $\mathcal{M}$ (resp. $\mathcal{M}^s$) denote the moduli space of smooth (resp. stable) cubic surfaces. Similarly, let $\mathcal{S}$ (resp. $\mathcal{S}^s$) denote the moduli space of smooth (resp. stable) $S_4$-symmetric cubic surfaces. Finally, let $\mathcal{H}^{S_4} \subset \CC\HH^1$ denote the ($S_4$-)symmetric discriminant locus of the period map.

\begin{theorem}\label{IntroMainThm1}
    There are analytic isomorphisms of orbifolds $\mathcal{S} \cong P\Gamma \backslash (\CC\HH^1 - \mathcal{H}^{S_4})$ and $\mathcal{S}^s \cong P\Gamma \backslash \mathbb{CH}^1$, where $\Gamma < \mathrm{U}(1,1)$ is an arithmetic lattice. Moreover the inclusion of moduli spaces $\mathcal{S} \to \mathcal{M}$ is compatible with the embedding of locally symmetric orbifolds $P\Gamma \backslash\mathbb{CH}^1 \to P\hat{\Gamma} \backslash\mathbb{CH}^4 \cong \mathcal{M}^s$.
\end{theorem}

For the precise statement of \autoref{IntroMainThm1}, its semistable extension, and its proof, see \autoref{MainTheoremActually}. Roughly, the idea behind the proof is to record the $S_4$-action in the period data and use the $S_4$-invariant subspace to define the period domain associated to symmetric cubic surfaces. We also determine the arithmetic group $\Gamma$ explicitly in \autoref{S4normArithmetic}. 

Let $\widetilde{\mathcal{M}}$ (resp. $\widetilde{\mathcal{S}}$) denote the space of (resp. symmetric) cubic surfaces equipped with a line. Recall that Jordan showed that the connected 27 lines cover $\widetilde{\mathcal{M}} \to \mathcal{M}$ has Galois group $W(E_6)$, the Weyl group of $E_6$.
Allcock--Carlson--Toledo recovered this fact Hodge-theoretically by considering an appropriate congruence cover of their uniformized moduli space $P \hat{\Gamma}\backslash \mathbb{CH}^4$ and using the exceptional isomorphism $W(E_6) \cong \PO(4,1,\FF_3)$. The following monodromy group result is an equivariant analog of Jordan's theorem for the symmetric $27$ lines cover $\widetilde{\mathcal{S}} \to \mathcal{S}$:

\begin{theorem}\label{IntroMainThm2}
    The (disconnected) symmetric $27$ lines cover $\widetilde{\mathcal{S}} \to \mathcal{S}$ of moduli spaces has monodromy group isomorphic to $S_4 \times K_4 < W(E_6)$. 
\end{theorem}

In the classical algebraic geometry, deformations and moduli problems were often studied using families of varieties with a chosen projective embedding; we refer to these universal families as \textit{parameter spaces}, from which the moduli spaces we are interested in can be obtained as a quotient by a projective linear action. For cubic surfaces, the monodromy groups of the 27 lines cover are identical for both the parameter and moduli spaces. In the symmetric locus, we prove that these two monodromy problems diverge.

Let $\mathcal{X}$ denote the moduli space of anti-canonically embedded smooth cubic surfaces and $\mathcal{Y}$ the moduli space of anti-canonically embedded symmetric cubic surfaces (see \autoref{sec:moduli} for more detail on how to build these spaces). Stack quotients, in the orbifold sense, of these spaces (by $\PGL_4$ and a group $N$ respectively, see \autoref{sec:moduli} for details) yield the moduli spaces $\mathcal{M}$ and $\mathcal{S}$. Each space admits a 27 lines cover ($\widetilde{\mathcal{X}}$ and $\widetilde{\mathcal{Y}}$, respectively) which each fits into a commutative diagram with deck groups given like so:
\begin{center}
    \begin{tikzcd}
\widetilde{\mathcal{X}} \arrow[d, "W(E_6)"'] \arrow[r, "\PGL_4\!\!\quot "] & \widetilde{\mathcal{M}} \arrow[d, "{\mathrm{PO}(4,1,\FF_3)}"] &  & \widetilde{\mathcal{Y}} \arrow[r, "N \quot"] \arrow[d, "M"'] & \widetilde{\mathcal{S}} \arrow[d, "S_4 \times K_4"] \\
\mathcal{X} \arrow[r, "\PGL_4 \!\!\quot"']                                 & \mathcal{M}                                                   &  & \mathcal{Y} \arrow[r, "N \quot "']                            & \mathcal{S}                                        
\end{tikzcd}
\end{center}

We determine the monodromy group $M$ which appears in the diagram above:

\begin{theorem}\label{thm:IntroMainThm3}
    The (disconnected) symmetric $27$ lines cover of parameter spaces $\widetilde{\mathcal{Y}} \to \mathcal{Y}$  has monodromy group $M$ isomorphic to centralizer $C_{W(E_6)}(S_4)$, which is a Klein 4-group.
\end{theorem}

The action of $K_4$ on $27$ labeled lines is explicitly worked out as permutations in \autoref{data:gens-K4-x-K4}. This allows us to completely characterize the covering space $\widetilde{\mathcal{Y}}$ --- it has 12 connected components, with each one corresponding to an explicit $K_4$-set; see \autoref{cor:incidence-var-structure} for details.

There are a few reasons why \autoref{thm:IntroMainThm3} is interesting. First off, the symmetric group $S_4$ and symmetric monodromy group $K_4$, thought of as subgroups of $W(E_6)$, intersect trivially --- this means that if we want to witness the $S_4$-action on a symmetric cubic surface through monodromy, we must leave the symmetric locus in the total parameter space. Second, the restrictions coming from Hodge theory constrain the monodromy group to a group of order 96 (this is \autoref{IntroMainThm2}). However, these restrictions provably do not suffice, as we can name explicit elements of this restricted subgroup that cannot arise via symmetric monodromy in the parameter space. This stands in direct contrast with reasoning used when studying similar problems, such as in \cite[Section 8]{RealACT}. A. Landesman pointed us to the root of this issue: the fibers of our stack quotient used to build the moduli space of symmetric cubic surfaces are not connected.

Finally, we give an explicit formula for the 27 lines on a symmetric cubic surface using only two square roots:

\begin{theorem}
    Given a general $S_4$-equivariant cubic surface $a \sum x_i^3 + b \sum x_i^2 x_j + c \sum x_i x_j x_k$ defined over $\mathbb{Q}$, its lines are all defined over the Klein four Galois extension $K = \mathbb{Q}(\sqrt{\sqrtone}, \sqrt{\sqrttwo})$, where
    \begin{align*}
        \sqrtone &= -{\left(9  a^{3} + 9  a^{2} b - 9  a b^{2} + 7  b^{3} - 3  a^{2} c - 6  a b c - 3  b^{2} c + 4  a c^{2}\right)} {\left(3  a + b - c\right)}\\
        \sqrttwo&=-{\left(3  a + b - c\right)} {\left(a + 3  b + c\right)}.
    \end{align*}
    Explicit formulas for lines in each orbit are given parametrically (as images of $\P^1$ with coordinates $[s:t]$) as follows:
    \begin{enumerate}
        \item[($S_4/C_2^o)$] The lines in this orbit are all the $S_4$-orbits of
        \begin{align*}
            \left[s + t : -s+t : \frac{9a^2-b^2-(3a-b)c + \sqrt{\sqrtone}}{6ab+2b^2 - 3(a+b)c+c^2}t : \frac{9a^2-b^2-(3a-b)c - \sqrt{\sqrtone}}{6ab+2b^2 - 3(a+b)c+c^2}t\right]
        \end{align*}
        \item[$(S_4/C_2^e)$] The lines in this orbit are all the $S_4$-orbits of
        \tiny
        \begin{align*}
            &\big[\frac{a-b-c+ \sqrt{\sqrttwo}}{2(a+b)}s + t : \frac{a-b-c+ \sqrt{\sqrttwo}}{2(a+b)}s - t :\\
            & s + \frac{(9a^2 - 6ab - b^2 + 2(3a + b)c - c^2) + (3a - 3b + c) \sqrt{\sqrttwo}}{2\sqrt{\sqrtone}} t : s - \frac{(9a^2 - 6ab - b^2 + 2(3a + b)c - c^2) - (3a - 3b + c) \sqrt{\sqrttwo}}{2\sqrt{\sqrtone}} t\big],
        \end{align*}
        \normalsize
        \item[$(S_4/D_8)$] The lines in this orbit are all the $S_4$-orbits of
        \[
            [s: -s : t : -t].
        \]
    \end{enumerate}
\end{theorem}
Given these symbolic formulas, we can now plot symmetric cubic surfaces together with their lines. A visual representation in \texttt{three.js} is available here:
\begin{center}
    \url{https://tbrazel.github.io/supplementary/eeg/s4_symmetric_cubics.html}
\end{center}

\subsection{Paper structure}
In \autoref{sec:moduli}, we review the construction of the (marked) moduli space of cubic surfaces, and explicitly realize the moduli space of symmetric cubic surfaces as a GIT quotient.
In \autoref{sec:ACT} we give an overview of the work of Allcock--Carlson--Toledo, including fundamental facts about framed cubic surfaces and their associated period maps, which yields their main theorem, a uniformization of the moduli of cubic surfaces by complex hyperbolic $4$-space.
In \autoref{sec:uniform}, we analogously define a period map for symmetric cubic surfaces via symmetric framings of cubic surfaces, and uniformize the moduli space of symmetric cubic surfaces by the complex hyperbolic line. This section finishes with a Hodge theoretic restriction of symmetric monodromy in the parameter space to a group of order $96$, which it turns out will properly contain the symmetric monodromy group.
In \autoref{sec:monodromy-group} we determine that the symmetric monodromy group is the Klein $4$-group $K_4$, and we describe how it acts on the $27$ lines. Moreover, we show that the $27$ lines cover over the symmetric locus $\widetilde{\mathcal{Y}} \to \mathcal{Y}$ splits into $12$ connected components, and explicitly determine the $K_4$-set structure associated to each connected component. \autoref{sec:repTheory} describes how to alternatively witness the symmetric monodromy group $K_4$ in $W(E_6)$ and $\PO(4,1,\FF_3)$ via representation theoretic constructions. Finally in \Cref{sec:formulas} we work out explicit formulas in radicals for lines on a symmetric cubic surface, leading to an elementary reproof of \cite[1.3]{EEG}, that a real smooth symmetric cubic surface contains only 3 or 27 real lines.

\autoref{sec:Appendix} contains explicit data regarding how the symmetric group $S_4$ and the generators of the symmetric monodromy group $K_4$ act on the $27$ lines on the Fermat cubic surface.



\subsection{Acknowledgements} We thank Benson Farb, Frank Sottile, and Jesse Wolfson for their interest and all independently asking us what the monodromy group is for the problem of finding lines on symmetric cubic surfaces. We thank Taylor Brysiewicz and  Kisun Lee for their help with homotopy continuation and certified tracking software. Thank you to Joe Harris, Aaron Landesman, Alberto Landi, Daniel Minahan, and Nick Salter for helpful conversations about this work. Lastly, we thank Daniel Allcock for his interest and elucidating answers regarding various group theoretic aspects of his work with James A. Carlson and Domingo Toledo. TB is supported by NSF Grant DMS-2303242, and SR is supported by NSF Grant DMS-2503485.

\subsection{Notation}
We use \texttt{\textbackslash{mathcal}} letters to indicate \textit{parameter spaces}, being both vector spaces parametrizing polynomials and moduli spaces of their vanishing loci.
\begin{center}
    \begin{tabular}{l | l }
    \hline
    \textbf{notation} & \textbf{meaning} \\
    \hline
    $\mathcal{W}$ & $\C[x_0, \ldots, x_3]_{(3)}$ \\
    $\mathcal{V}$ & $\C[x_0, \ldots, x_3]_{(3)}^{S_4}$ \\
    $\mathcal{X}$ & moduli of anti-canonically embedded cubic surfaces \\
    $\mathcal{Y}$ & moduli of anti-canonically embedded symmetric cubic surfaces \\
    $\mathcal{M}$ & moduli of cubic surfaces \\
    $\mathcal{S}$ & moduli of symmetric cubic surfaces \\
    \hline
    \textbf{decoration} & \textbf{meaning} \\
    \hline
    $(-)^{\sm}$ or no decoration & moduli of smooth objects \\
    $(-)^{s}$ & moduli of stable objects \\
    $(-)^{ss}$ & moduli of semistable objects \\
    $\widetilde{(-)}$ & moduli of cubic surfaces equipped with a line \\
    $\widehat{(-)}$ & moduli of marked cubic surfaces
    \end{tabular}
\end{center}

\section{Moduli constructions}\label{sec:moduli}

The content of this section is standard and well-known --- see \cite{ZhengOrbifold} for example. We will review how to construct the moduli space of (marked) smooth cubic surfaces as a GIT quotient, and then analogously construct the moduli space of (marked) symmetric cubic surfaces. We end this section with a quick discussion of stability and semistability of cubic surfaces, concluding with the facts that the Cayley nodal cubic is the only non-smooth symmetric stable cubic surface, and that the tricuspidal cubic is the only non-stable symmetric semistable cubic surface.

\subsection{Parameter space of cubic surfaces}
Let $\mathcal{W}=\CC[z_0,z_1,z_2, z_{3}]_{(3)}$ denote the $20$-dimensional vector space of degree $3$ homogeneous polynomials in $4$ variables. Every $f \in \mathcal{W}\backslash \{0\}$ defines a cubic surface $Z(f)$ in $\PP^{3}$, and two elements $f_1$ and $f_2$ determine the same surface $Z(f_1) = Z(f_2)$ if and only if $f_1 = \lambda f_2$ for some $\lambda \in \CC^*$. Thus $\PP (\mathcal{W}) \cong \PP^{19}$ can be naturally thought of as the parameter space of cubic surfaces in $\PP^{3}$. 

There is a linear action of $\SL(4,\CC)$ on $\mathcal{W}$ given by $g \cdot f := f \circ g^{-1}$. This induces a left action of $\PGL(4,\CC)$ on the projectivization $\PP(\mathcal{W})$.

\begin{definition} Consider the left $\SL(4,\C)$-action on $\mathcal{W}$ induced by permuting coordinates on $\PP^3$. For $f\in \mathcal{W}$, we say that $f$ is
\begin{enumerate}
    \item \textit{smooth} if its associated cubic surface is smooth,
    \item \textit{stable} if the orbit $\SL(4,\C)\cdot f$ is closed, and the stabilizer subgroup is finite,
    \item \textit{semistable} if $0$ is not in the closure of the orbit $\SL(4,\C)\cdot f$.
\end{enumerate}
\end{definition}

We denote by $\mathcal{W}^{\sm}$ (respectively $\mathcal{W}^{s}$, and $\mathcal{W}^{ss}$) the subsets of $\mathcal{W}$ corresponding to smooth cubic surfaces (respectively stable, and semistable). It is classically known that we have containments
\begin{align*}
    \mathcal{W}^{\sm}\subseteq \mathcal{W}^{s}\subseteq \mathcal{W}^{ss}.
\end{align*}
The action of $\SL(4,\CC)$ on each of these loci descends to an action of $\PGL(4,\CC)$ on their projectivizations. We can take the respective GIT quotients to construct various moduli spaces of cubic surfaces.

\begin{definition} We denote by
\begin{align*}
    \mathcal{M}^{\sm} &:=\PGL(4,\CC)\quot \PP(\mathcal{W}^{\sm}) ,\\
    \mathcal{M}^{s} &:=\PGL(4,\CC)\quot \PP(\mathcal{W}^{s}) ,\\
    \mathcal{M}^{ss} &:=\PGL(4,\CC)\quot \PP(\mathcal{W}^{ss}) ,
\end{align*}
the \textit{moduli space of smooth/stable/semistable cubic surfaces}.
\end{definition}

\begin{convention} When we write a moduli space without a superscript, e.g. $\mathcal{M}$, we implicitly mean the moduli of smooth objects.
\end{convention}

The (finite) pointwise stabilizer subgroups in $\PGL(4,\CC)$ of cubic forms yield automorphisms of the corresponding cubic surface, and thus naturally gives us an orbifold structure on the GIT quotient $\mathcal{M}$. The same will hold true in the symmetric locus, and so we will regard those quotients as orbifolds as well.

The following classical result characterizes stable and semistable cubic surfaces by their singularities.

\begin{theorem}[{Hilbert, \cite{hilbertVollen}}] 
A cubic surface is stable if and only if it its singularities are ordinary nodes. A cubic surface is semi-stable if and only if its singularities are ordinary nodes or $A_2$ singularities.
\end{theorem}

\begin{lemma}[{\cite[4.6]{AllcockCarlsonToledo}}]\label{lem:tricuspidal}
The cubic form $z_0^3 - z_1 z_2 z_3$ is the unique closed $\SL(4,\C)$-orbit of semistable non-stable cubic surfaces.
\end{lemma}

Since points in the GIT quotient $\mathcal{M}^{ss}$ correspond to closed orbits, this indicates that there is a unique point in the moduli space of semistable cubic surfaces corresponding to a point which is not stable. This is given by the unique \textit{tricuspidal cubic surface}, defined by the equation mentioned, and pictured in \autoref{fig:cayley}.

Let $\widetilde{\mathcal{W}^{\sm}}$ denote the parameter space of smooth cubic forms equipped with an incident line. Concretely, this is the incidence variety
$$\widetilde{\mathcal{W}^{\sm}} = \{(f,\ell) \in \mathcal{W}^{\sm} \times \Gr(2,4) : \ell \subset Z(f)\}.$$
The GIT quotient $\widetilde{\mathcal{M}} = \PGL(4,\CC)\backslash\!\!\backslash \PP \widetilde{\mathcal{W}^{\sm}}$ is the moduli space of smooth cubic surfaces equipped with a line. Since the natural projection $\widetilde{\mathcal{W}^{\sm}} \to \mathcal{W}^{\sm}$ is $\PGL(4,\CC)$-equivariant, it descends to a map of moduli spaces $\widetilde{\mathcal{M}} \to \mathcal{M}$.

\subsection{Marked parameter space of cubic surfaces}

Recall that a free finitely generated $\ZZ$-module $L$ equipped with an integral symmetric (resp. symplectic) non-degenerate bilinear form $q$ defines a \textit{symmetric (resp. symplectic) lattice structure} $(L,q)$. The lattice structure on the intersection form of a smooth cubic surface is classically obtained by viewing the surface as a blowup of the projective plane at six points.

\begin{proposition}
Let $X = V(f) \subset \PP^3$ denote a smooth cubic surface determined by some cubic form $f \in \mathcal{W}^{\sm}$. Then $H= H^2 (X,\ZZ)$ is a free $\ZZ$-module of rank $7$, and the cup product $\langle \cdot,\cdot\rangle$ determines a signature $(1,6)$ symmetric unimodular lattice structure $(H,\langle \cdot,\cdot\rangle)$.
\end{proposition}

Let $\eta_X \in H$ denote the canonical class on $X$ and $(L,q) \cong \langle 1 \rangle \oplus \langle -1 \rangle^{\oplus 6}$ be an abstract lattice isomorphic to $(H ,\langle \cdot,\cdot\rangle)$. Fix some $\eta \in L$ so that $(L,q,\eta) \cong (H,\langle \cdot,\cdot\rangle,\eta_X)$.

\begin{definition}\label{MarkingsAndMarkedModuli}
    A \textit{marking} of a smooth cubic surface $X$ is an isomorphism of lattices
    $$\phi : (H, \langle \cdot,\cdot\rangle, \eta_X) \to (L,q, \eta).$$
   We say a cubic form with marking $(f_1,\phi_1)$ is equivalent to the pair $(f_2,\phi_2)$ if there exists some $g \in \PGL(4,\CC)$ so that $g(f_1) = f_2$ and $\phi_2 = \phi_1 \circ g^*$. We will let $\widehat{\mathcal{W}^{\sm}}$ denote the parameter space of marked smooth cubic forms, which is naturally a complex manifold \cite[3.2]{AllcockCarlsonToledo}.
\end{definition}

Let $\widehat{\mathcal{M}}$ denote the GIT quotient $\PGL(4,\CC)\backslash\!\!\backslash \PP\widehat{\mathcal{W}^{\sm}}$. We refer to this as the \textit{moduli space of smooth marked cubic surfaces}. As cubic surfaces vary, their markings will vary as well. Since any two markings differ by an automorphism of the abstract lattice $(L,q,\eta)$, we obtain a representation of the fundamental group of the moduli space of smooth marked cubic surfaces. The following is a relevant consequence of work of Beauville on monodromy in the universal family of degree $d$ hypersurfaces which was classically known for cubic surfaces \cite{BeauvilleMonodromy}:

\begin{proposition}\label{WholeBeauvilleMonodromy}
    The space $\widehat{\mathcal{M}}$ is a connected, Hausdorff complex manifold which is a covering space of $\mathcal{M}$. Moreover, the monodromy representation
    $$\pi_1(\mathcal{M},X) \to \Aut(H^2 (X,\ZZ),\eta_X)$$
    is surjective and has image isomorphic to the Weyl group of the root lattice $E_6$, denoted $W(E_6)$.
\end{proposition}

It is also classically known that the moduli space of marked cubic surfaces $\widehat{\mathcal{M}}$ is isomorphic to the moduli space of cubic surfaces equipped with six ordered skew lines \cite[pg. 19]{BeauvilleBourbaki}. We shall freely identify these spaces. Moreover, the cover $\widehat{\mathcal{M}} \to \mathcal{M}$ is the normal closure of the 27 lines cover $\widetilde{\mathcal{M}} \to \mathcal{M}$; we will work on these marked moduli spaces (and their symmetric analogs) to determine our desired Galois groups.

\subsection{Parameter space of symmetric cubic surfaces}

Recall that a degree $d$ homogeneous polynomial $f(z_0,\dots,z_{n})$ is \textit{symmetric} if it is invariant under natural $S_{n+1}$-action on $z_0,\dots,z_{n}$, i.e. $f(z_0 ,\dots, z_n) = f(z_{\sigma(0)},\dots,z_{\sigma(n)})$ for all $\sigma \in S_{n+1}$. When $n >d$, the vector space $\CC[z_0 ,\dots,z_{n}]_{(d)}^{S_{n+1}}$ of symmetric homogeneous degree $d$ polynomials in $n+1$ variables is $p(d)$-dimensional, where $p(d)$ denotes the number of partitions of $d$. A basis will be denoted by $\{m_\alpha\}$, where $m_\alpha$ is a homogeneous symmetric polynomial indexed by the partitions $\alpha\vdash d$.

In the case of symmetric cubic forms in $4$ variables, the vector space $\mathcal{V} := \mathcal{W}^{S_4}$ admits a basis of the form
\begin{align*}
    m_{3}(z_0,z_1 ,z_2, z_3) &= \sum z_i^3,\\
    m_{21}(z_0,z_1 ,z_2, z_3) &= \sum z_i^2 z_j,\\
    m_{111}(z_0,z_1 ,z_2, z_3) &= \sum z_i z_jz_k,
\end{align*}
so any symmetric cubic form $f$ in $4$ variables can be uniquely written as a linear combination
$$f =  a \cdot m_{3} + b\cdot m_{21} + c \cdot m_{111}.$$
We see that the the parameter space of symmetric cubic forms $\PP (\mathcal{V}) =\PP^2$ embeds linearly into the parameter space of cubic forms $\PP(\mathcal{W})$.
Define $ \Delta^{S_4}$ to be the symmetric discriminant curve, which corresponds to the locus of $S_4$-invariant singular cubic forms in the parameter space $\mathcal{Y}$.

In order to form a GIT quotient parametrizing a moduli space of symmetric cubic surfaces, we need to understand how the action of $\PGL(4,\CC)$ preserves or fails to preserve the symmetry of the associated cubic surface. The following is a basic algebra fact that will be relevant to much of what follows:

\begin{proposition}\label{prop:S4normalizercentral}
    Let $S_4$ be a subgroup of any group $G$. Then the normalizer $N_G(S_4)$ is generated by $S_4$ and its centralizer $C_G(S_4)$.\footnote{The general fact we are using is that every automorphism of a complete group is inner. Thus a complete subgroup of any group has normalizer generated by the subgroup and its centralizer. The argument we give works \textit{mutatis mutandis}.}
\end{proposition}
\begin{proof}
    Given any $g \in N_G(S_4)$, we have $g \sigma g^{-1} \in S_4$ for all $\sigma \in S_4$. Thus conjugation by $g$ defines an automorphism of $S_4$. Recall that $S_n$ is a complete group for $n \neq 2,6$, and so every automorphism of $S_4$ is an inner automorphism. This implies that for each $g \in N_G(S_4)$, there exists some $\eta \in S_4$ such that
    $$g \sigma g^{-1} = \eta \sigma \eta^{-1} \Leftrightarrow \sigma = \eta^{-1}g \sigma g^{-1}\eta = \eta^{-1}g \sigma (\eta^{-1}g)^{-1} $$
    for all $\sigma \in S_4$. Thus $\eta^{-1} g \in C_G(S_4)$, and so $g \in C_G(S_4) \cdot S_4$. This proves the claim.
\end{proof}

\begin{proposition}\label{PGLnormalizerofS4}
    The normalizer of the permutation subgroup $S_4$ in $\PGL(4,\CC)$ is 
    $$N_{\PGL(4,\CC)} (S_4) \cong \left\{
    {\scalebox{0.75}{$\begin{pmatrix*}
        \lambda & 1 & 1 & 1\\\
        1 & \lambda & 1 & 1\\
        1 & 1 & \lambda & 1\\
        1 & 1 & 1 & \lambda
    \end{pmatrix*}$}} 
    \cdot P : \lambda\notin\{1,-3\},\, P \in S_4 \leq \PGL(4,\CC) \right\}$$
\end{proposition}
\begin{proof}
    By \autoref{prop:S4normalizercentral}, it suffices to determine the centralizer of $S_4$ in $\PGL(4,\CC)$. One can then calculate that the subgroup of permutation matrices in $\GL(4,\CC)$ is centralized by matrices of the form
    $$C(a,b)= {\scalebox{0.75}{$\begin{pmatrix*}
        a & b & b & b\\
        b & a & b & b\\
        b & b & a & b\\
        b & b & b & a
    \end{pmatrix*}$}}$$
    where $ b \neq a, a \neq -3 b$.

    Let $\varphi : \GL(4,\CC) \to \PGL(4,\CC)$ be the projectivization homomorphism. Since the permutation matrices intersect the center of $\GL(4,\CC)$ trivially, we have $\varphi(S_4) \cong S_4$. To determine the rest of the image of $N_{\PGL(4,\CC)}(S_4)$, we break into two cases, when $b = 0$ or $b \neq 0$. If $b = 0$ then the matrices $C(a,0)$ are scalar and form the kernel of $\varphi$. If $b \neq 0$, then $\varphi(C(a,b)) = \varphi(C(a/b,1))$. For all $a \in \CC \backslash \{1,-3\}$, the matrix $C(a/b,1)$ induces a nontrivial automorphism of $\PP^3$ and thus does not lie in the kernel of $\varphi$. Thus by letting $\lambda = a/b$, we obtain that the normalizer of $S_4 < \PGL(4,\CC)$ is the subgroup stated in the proposition.
\end{proof}

This allows us to define the moduli of smooth symmetric cubic surfaces.

\begin{proposition}\label{prop:S-definition}
The GIT quotient $\mathcal{S} = N_{\PGL(4,\CC)} (S_4) \quot \mathcal{Y}^{\sm}$ is the moduli space of smooth symmetric cubic surfaces.
\end{proposition}
\begin{proof}
    Suppose that two symmetric cubic forms $f_1,f_2 \in \mathcal{Y}^{\sm}$ determine isomorphic symmetric cubic surfaces $X = Z(f_1)$ and $Y = Z(f_2)$. Since such an isomorphism $\varphi: X \xrightarrow[]{\sim} Y$ of varieties preserves their respective canonical classes $K_X$ and $K_Y$, the map extends to respect their anticanonical embeddings into $\PP^3$. Thus such an isomorphism $\varphi$ must be the restriction of a linear automorphism coming from the ambient projective space $\PP^3$. Moreover, the automorphism groups of $X$ and $Y$ must be preserved under such an isomorphism, and so the $S_4$-action on $X$ must be sent to the $S_4$-action on $Y$. Thus two symmetric cubic surfaces are projectively equivalent when their symmetric cubic forms differ by an element of the normalizer of $S_4 < \PGL(4,\CC)$, which was explicitly calculated in \autoref{PGLnormalizerofS4}.
\end{proof}

\subsection{Symmetry and stability} Having defined the moduli space of smooth symmetric cubic surfaces $\mathcal{S}$ in \autoref{prop:S-definition}, we would like to define the analogous moduli spaces of stable and semistable symmetric cubic surfaces. In order to do this, we first must explore how (semi)stability interacts with symmetry.

\begin{proposition}\label{prop:singularities-on-semistable-singular-cubic}
Let $f \in \mathcal{V}$ be a nonzero symmetric homogeneous form defining a semistable cubic surface. Then the singularities of $V(f)$ are either $4A_1$ or $3A_2$.
\end{proposition}
\begin{proof} We first see that if $X = V(f)$, then the geometric $S_4$ action it inherits by symmetry is actually a subgroup of the automorphism group. This is clear if $X$ is smooth, since the cubic surface is anticanonically embedded, but a small argument is needed if $X$ isn't smooth. Suppose towards a contradiction that any non-trivial element $g\in S_4$ acted trivially on $X$. Then $X$ would lie in the $g$-fixed subspace of $\P^3$, which is a hyperplane or intersection of hyperplanes, a contradiction.

Since a semistable cubic surface is normal by Serre's criterion \cite[5.10]{grothendieck1965elements}, we can refer to the classification of automorphism groups of normal cubic surfaces due to Sakamaki \cite[Table~3]{Sakamaki}. It is clear that $S_4$ cannot be a subgroup of any of the automorphism groups except $4A_1$ where it is equality, and $3A_2$, where we make use of the semidirect product $K_4\rtimes S_3 \cong S_4$.
\end{proof}

\begin{figure}[h]
\includegraphics[width=0.4\linewidth]{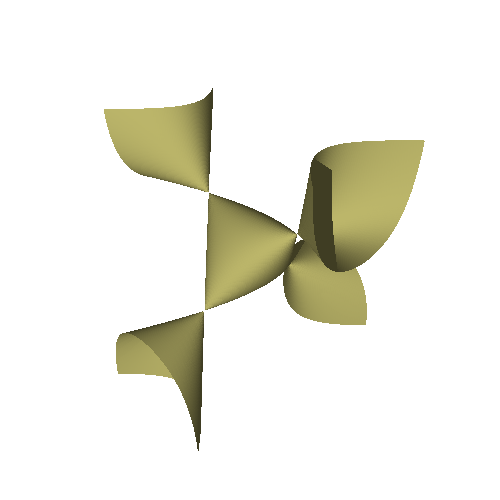}
\includegraphics[width=0.4\linewidth]{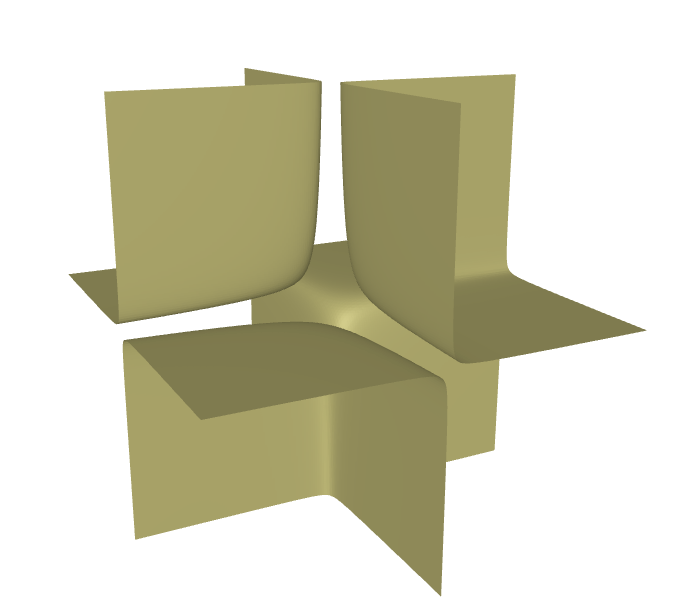}
  \centering
  \caption{Left: the Cayley nodal cubic surface. Right: the tricuspidal cubic surface}
  \label{fig:cayley}
\end{figure}

We now look to see if any such symmetric singular cubic surfaces do exist. One of the most famous singular cubic surfaces is symmetric:

\begin{definition} The \textit{Cayley nodal cubic surface}, defined by the elementary symmetric homogeneous form $m_{111}$ is a singular cubic surface with four nodes. Its automorphism group is $S_4$, which permutes these four nodes \cite{Sakamaki}. It is pictured in \autoref{fig:cayley}.
\end{definition}

Conveniently, the normalizer of $S_4$ in $\PGL(4,\CC)$ appears in the following proposition, which characterizes the Cayley cubic surface as the unique cubic surface with four nodes (c.f. \cite{BruceWall}).

\begin{proposition} Let $f\in \mathcal{W}$ be a nonzero form defining a cubic surface with four nodes. Then there exists a \textit{unique} change of coordinates $g\in N_{\PGL(4,\CC)}(S_4)$ so that $g\cdot f$ is the Cayley nodal cubic.
\end{proposition}
\begin{proof} Given any other cubic surface with four nodes, there is a projective change of coordinates turning it into the Cayley nodal cubic by sending the four nodes to the four nodes of the Cayley cubic. This change of coordinates is unique up to permutation of the nodes since $\PGL(4,\CC)$ is simply 4-transitive. However there is a unique automorphism of the Cayley nodal cubic corresponding to any permutation of the nodes, since its automorphism group is the symmetric group $S_4$.
\end{proof}

This has an immediate corollary to our study of stable symmetric cubic surfaces, which allows us to understand the moduli space.
\begin{corollary}
Any symmetric stable cubic surface which is not smooth lies in the $N_{\PGL(4,\CC)}(S_4)$ orbit of the Cayley nodal cubic surface.
\end{corollary}

\begin{corollary}\label{TheGITConstructionOfStableSymmetricCubics}
The moduli space of stable symmetric cubic surfaces
\begin{align*}
    \mathcal{S}^s = N_{\PGL(4,\CC)} (S_4)  \backslash\!\!\backslash \PP(\mathcal{V}^s)
\end{align*}
has the property that $\mathcal{S} \subseteq\mathcal{S}^s$, and $\mathcal{S}^s \minus\mathcal{S} = \{C\}$ is one point, which is the Cayley nodal cubic.
\end{corollary}

What about semistability? By \autoref{prop:singularities-on-semistable-singular-cubic} any semistable symmetric cubic which isn't stable must have three cusps, and we know there is a unique semistable non-stable cubic surface in the moduli space $\mathcal{M}^{ss}$ by \autoref{lem:tricuspidal}. So it suffices to check if there is any projective change of coordinates exhibiting the tricuspidal cubic surface as a symmetric homogeneous form.

\begin{computation} The tricuspidal cubic surface $z_0^3 - z_1 z_2 z_3$ is projectively equivalent to the homogeneous form $4m_{21} + 4m_{111}$ via the change of basis matrix
\begin{align*}
    \begin{pmatrix}
        1 & 1 & 1 & 1 \\
        1 & 1 & -1 & -1 \\
        1 & -1 & 1 & -1 \\
        1 & -1 & -1 & 1.
    \end{pmatrix}
\end{align*}
\end{computation}
\begin{proof}[Why this works] Since $S_4$ acts faithfully on the cubic surface, it must map singularities to singularities, hence if a symmetric cubic surface has three cusps, they form an $S_4$-set. Viewing $\P^3$ as an $S_4$-space, we see there is a unique point in $\P^3$ with full isotropy $S_4$, and no points with isotropy $A_4$. Hence the three cusps must form a transitive $S_4$-set, isomorphic to $S_4/D_8$. We check that there is a unique such collection of three points in $\P^3$, namely $[1:1:-1:-1]$, $[1:-1:1:-1]$, and $[1:-1:-1:1]$. Since $\PGL(4,\CC)$ is 3-transitive, if we can the tricuspidal cubic above into a symmetric form, we must map its cusps $[0:1:0:0]$, $[0:0:1:0]$, and $[0:0:0:1]$ to the points forming the $S_4/D_8$ orbit above, hence the three rightmost columns in the matrix we found. A computation then forces the first column to consist of all 1's.
\end{proof}

\begin{corollary} The moduli space of semistable symmetric cubic surfaces
\begin{align*}
    \mathcal{S}^{ss} := N_{\PGL(4,\CC)}(S_4) \quot \PP(\mathcal{V}^{ss})
\end{align*}
has exactly one point not in the stable moduli space, corresponding to the tricuspidal curve $m_{21} + m_{111}$.
\end{corollary}

The inclusion $\mathcal{V} \to \mathcal{W}$ of parameter spaces descends to a inclusion $\iota : \mathcal{S} \to \mathcal{M}$ of the moduli of symmetric cubic surfaces into the whole moduli space. The pullback $\iota^* \widehat{\mathcal{M}}$ is then the moduli space of smooth marked symmetric cubic surfaces $\widehat{\mathcal{S}}$. Then just as before, the pullback construction let's us conclude that $\iota^* \widehat{\mathcal{M}}$ is isomorphic to the moduli space of symmetric cubic surfaces equipped with six ordered skew lines, which we will also denote by $\widehat{\mathcal{S}}$.

\section{Reviewing Allcock--Carlson--Toledo}\label{sec:ACT}

In this section, we outline the construction of the Allcock--Carlson--Toledo period map \cite{AllcockCarlsonToledo}.
Let $S=Z(f)$ denote a cubic surface in $\PP^3$. Since $H^2(S,\ZZ)$ admits a type $(1, 1)$ Hodge structure, the natural period map is constant. In their seminal paper, Allcock--Carlson--Toledo showed that there is an weight 3 Hodge structure associated to $S$ whose periods entirely capture its geometry. This is the Hodge structure of the cyclic cubic threefold $T$, realized as a degree $3$ cover of $\PP^3$ with branch locus $S$. In coordinates,
$$T = \{t^3 = f(z_0,z_1,z_2,z_3)\} \subset \PP^4,$$
and the deck group $\langle \tau  \rangle$ of the cover acts on $T$ by multiplying the $t$-coordinate by the 3rd root of unity $\omega$. The pair $(H^3(T,\ZZ),\tau)$ forms a so-called \textit{Eisenstein Hodge structure}. The period domain for such Hodge
structures is complex hyperbolic $4$-space $\CC\HH^4$. There is a natural period
map from the moduli of smooth cubic surfaces $\mathcal{M}$ to an arithmetic quotient of $\CC\HH^4$ by $P\hat{\Gamma} = \mathrm{PU}(4,1,\ZZ[\omega])$:
$$\mathcal{P} : \mathcal{M} \to P\hat{\Gamma}\backslash\CC\HH^4.$$
By the Riemann extension theorem, the map $\mathcal{P}$ extends uniquely to the stable cubic surfaces $\mathcal{M}^s$. The main theorem of \cite{AllcockCarlsonToledo} is that $\mathcal{P}$ is a biholomorphism of analytic spaces $\mathcal{M}^s \cong P\hat{\Gamma}\backslash\CC\HH^4$, and moreover it is an isomorphism of orbifolds. After reviewing their work, we will build an analogous uniformization of the moduli space of symmetric cubic surfaces. 

\subsection{Basic cohomology knowledge}
Given a smooth cubic surface $S$ defined by the cubic form $f$, we associate to it the cyclic cubic $3$-fold $T:=\{t^3 = f\}$ which defines a degree 3 branched covered $\PP^3$ over $S$. The Lefschetz hyperplane theorem and Poincar\'{e} duality tells us that $T$ and $\PP^3$ have the same cohomology away from the middle degree 3, with Hodge numbers equal to $h^{i,i}(T) = 1$ for $i=0,1,2,3$. An Euler characteristic calculation tells us that $H^3(T,\ZZ)$ is rank 10, and since there are no holomorphic 3-forms on $T$, i.e. $h^{3,0} (T) = h^{0,3} (T) = 0$, the middle Hodge numbers of $T$ are $h^{2,1} = h^{1,2} = 5$. 

\subsection{The module structure on cohomology}

The deck group $\langle \tau \rangle$ of the triple branched cover $T \to \PP^3$ acts on $T$ by multiplying the $t$-coordinate by a 3rd root of unity $\omega$. Since fixed vectors of the induced action on cohomology $H^3(T,\ZZ)$ must come from $H^3(\PP^3,\ZZ) = 0$ via transfer, $\tau$ acts on the (real) cohomology of $T$ without fixed points. Thus the minimal polynomial of the $\tau$-action on $H^3(T)$ is $z^2 + z + 1$, and $H^3(T,\ZZ)$ inherits the structure of a free $\ZZ[\omega]$-module of dimension five.

The Hodge decomposition on $H^3(T,\CC)$ forms a direct sum decomposition
$$H^3(T,\CC) = H^{2,1}(T) \oplus H^{1,2}(T),$$
where the two summands are isomorphic dimension five vector spaces and are exchanged by complex conjugation. Since $\tau$ acts holomorphically on $T$, $\omega$ acts on $H^3(T,\ZZ)$ as a real operator, and so it preserves the Hodge decomposition. Thus the eigenspace decomposition $H^3(T) = H_\omega^3(T) \oplus H_{\bar{\omega}}^3(T)$ is compatible with the Hodge decomposition. Selecting the $\bar{\omega}$-summand, we get a Hodge-eigenspace direct sum decomposition
$$H_{\bar{\omega}}^{3}(T) = H_{\bar{\omega}}^{2,1}(T) \oplus H_{{\bar{\omega}}}^{1,2}(T).$$
There is a naturally associated Hermitian $\ZZ[\omega]$-valued form $h$ on $H^3(T,\ZZ)$ coming from the cup product $\langle \cdot, \cdot \rangle$:
$$h(a,b) := \frac{\langle (\tau - \tau^{-1}) a,  b\rangle - (\omega - \omega^{-1}) \langle a, b\rangle}{2}.$$
With respect to this form $h$, the Hermitian pair $(H_{\bar{\omega}}^3 (T), h)$ is a signature $(4,1)$ complex inner product space. A key point in \cite{AllcockCarlsonToledo} is that this direct sum decomposition is orthogonal with respect to $h$, $H_{\bar{\omega}}^{2,1}(T)$ is $1$-dimensional and negative-definite with respect to $h$, and that $H_{\bar{\omega}}^{1,2}(T)$ is $4$-dimensional and positive-definite with respect to $h$.

\subsection{Moduli of framed cubic surfaces}

We have already defined markings and constructed the moduli space of marked cubic surfaces. A related but slightly different idea must be studied, which is the notion of a \textit{framing}; this is a marking of the cohomology of the cyclic cubic threefold $T$ associated to a cubic surface $S$.

\begin{definition}
    A \textit{framing} of the cubic surface $S$ is an isometry $\psi$ from the $\ZZ[\omega]$-lattice $(H^3(T,\ZZ),h)$ to the abstract indefinite $\ZZ[\omega]$-lattice $\Lambda = \ZZ[\omega]^{4,1}$ (recall that $\Lambda$ is unique up to isometry \cite{Allcock}). Two framings $(S_1, \psi_1)$ and $(S_2, \psi_2)$ are equivalent if there exists some $g \in \PGL(4,\CC)$ such that $\psi_1 = \tilde{g}^* \psi_2$. The space of equivalence classes of framed cubic surfaces $(S,\psi)$ is denoted by $\mathcal{F}$.
\end{definition}

Allcock--Carlson--Toledo showed that $\mathcal{F}$ is a complex manifold. There is a natural action of the arithmetic group $\hat{\Gamma} = \mathrm{U}(4,1,\ZZ[\omega])$ on the space of framed cubic surfaces by $\gamma \cdot (S, \psi) = (S, \psi \circ \gamma^{-1})$.

\subsection{The period mapping}
The space of negative lines in the space $\CC^{4,1} = \Lambda \otimes_{\ZZ[\omega]} \CC$ is the complex hyperbolic $4$-space $\CC\HH^4$. We can now define the period mapping of framed cubic surfaces $\widetilde{\mathcal{P}} : \mathcal{F} \to \CC\HH^4$ to be
$$\widetilde{\mathcal{P}}:(S,\psi) \longmapsto \PP(\psi(H_{\bar{\omega}}^{2,1}(T))) \in \CC\HH^4 \subset \CC\PP^{4}.$$
The main theorem of \cite{AllcockCarlsonToledo} is that this map is an open embedding, and descends through the $\hat{\Gamma}$-quotient to an isomorphism $\mathcal{P}$ of moduli spaces $\mathcal{M} \cong P\hat{\Gamma} \backslash (\CC\HH^4-\mathcal{H})$, where $\mathcal{H}$ is the locally finite hyperplane arrangement determined by reflections over short roots $\delta\in \Lambda$:
$$\mathcal{H} = \bigcup_{h(\delta,\delta) = 1} \mathrm{Fix}(\mathrm{Ref}_\delta).$$
This isomorphism extends to the moduli space of stable cubic surfaces $\mathcal{M}^s \cong P\hat{\Gamma} \backslash \CC\HH^4$. They remark that this map extends to an analytic isomorphism of Deligne--Mumford stacks; for more details on this stacky structure, see the papers of Kudla--Rapoport \cite{KudlaRapoport} and Zheng \cite{ZhengOrbifold}.

\begin{remark}
    Note that originally Allcock--Carlson--Toledo used the negative-definite line $H_\omega^{1,2}(T)$ to define the period data of cubic surfaces, but this makes the period mapping anti-holomorphic, as pointed out by Beauville \cite{BeauvilleBourbaki}. This is why we adopted the $\bar{\omega}$ convention instead.
\end{remark}

Finally, consider the group homomorphism $\ZZ[\omega] \to \mathbb{F}_3$ which sends $\omega$ to $1$. Then $\ZZ[\omega]^{4,1} \otimes_{\ZZ[\omega]}\FF_3 \cong \FF_3^{4,1}$, which induces a homomorphism 
$$\varphi : \hat{\Gamma} \to \mathrm{PO}(4,1,\FF_3).$$
Allcock--Carlson--Toledo showed that this homomorphism is surjective, with kernel denoted by $\hat{\Gamma}'$. The following is well-known which we include for the sake of completeness (see \cite[Section 2.12]{AllcockCarlsonToledo} and \cite[pg. 26]{TheATLAS}):

\begin{proposition}\label{ExceptionalIsomorphismWE6}
    There is an exceptional isomorphism of finite groups $W(E_6) \cong \mathrm{PO}(4,1,\FF_3)$. 
\end{proposition}
\begin{proof}[Proof sketch]
    The Weyl group $W(E_6)$ acts on the $6$-dimensional root lattice $E_6$. This has index $3$ in the weight lattice. The root lattice modulo $3$ times the weight lattice is $5$-dimensional vector space over $\FF_3$. This space inherits an inner product $q$ from the root lattice by reducing mod $3$ the inner product of lattice vectors. Thus every element of $W(E_6)$ descends to an automorphism of $\FF_3^5$ which preserves this non-degenerate symmetric bilinear form $q$. This yields the homomorphism $W(E_6) \to \mathrm{O}(q,\FF_3)$. Post-composing with the projectivization, we obtain the group homomorphism 
    $$W(E_6) \to \mathrm{PO}(q,\FF_3).$$
    Any two non-degenerate symmetric bilinear forms $q$ and $q'$ are equivalent over $\FF_3$ \cite{HusemollerMilnor}, and so the group $\mathrm{PO}(q,\FF_3)$ and $\mathrm{PO}(4,1,\FF_3)$ isomorphic. One can then calculate that $|W(E_6)| = |\mathrm{PO}(q,\FF_3)| = 51840$. By almost simplicity of $W(E_6)$ and non-triviality of this homomorphism, this map is an isomorphism.
\end{proof}

By \autoref{ExceptionalIsomorphismWE6}, we have that $\hat{\Gamma}/\hat{\Gamma}' \cong W(E_6)$, and so the Galois cover of $\mathcal{M}$ that $\hat{\Gamma}' < \hat{\Gamma}$ corresponds to is the space of cubic surfaces equipped with a line $\widehat{\mathcal{M}}$. Once we have analogously uniformized the moduli space of symmetric cubic surfaces, a special symmetric subgroup of $\mathrm{PO}(4,1,\FF_3)$ will be determined that is the monodromy group of the cover $\widetilde{\mathcal{S}} \to \mathcal{S}$.

\section{Uniformization of the symmetric moduli space}\label{sec:uniform}

As was reviewed in the previous section, Allcock--Carlson--Toledo showed that the moduli space of stable cubic surfaces $\mathcal{M}^{s}$ admits the structure of a complex ball quotient. Specializing to the stable symmetric locus $\mathcal{S}^{s}$, we will study this space through Hodge theory and analogously realize it as a ball quotient. Although their stated goals and some of the technology used are different, we take inspiration from and owe an intellectual debt to the work of Yu--Zheng \cite{YuZheng}.

To understand the moduli space of symmetric cubic surfaces, we need to understand how $S_4$ acts on the period data. More specifically, we want to understand the $\bar{\omega}$-eigenspace in cohomology $H_{\bar{\omega}}^{3}(T)$ as an $S_4$-representation. Any $S_4$-equivariant automorphism of the cubic surface $S$ will lift to give a nontrivial $S_4$-equivariant action on $T$, thus the cohomology group $H_{\bar{\omega}}^{3}(T)$ while preserving the Hodge decomposition. To understand this action, we will express the periods of a given symmetric cubic surface in terms of differential forms. 

\subsection{Residue calculus and symmetric cubic $3$-folds}

Griffiths' theory of residues \cite{GriffithsRational} will help us make explicit the action of $S_4$ on the invariant cohomology of a cyclic cubic $3$-fold. His foundational work on rational integrals allows us to assert the following:

\begin{proposition}
    Let $S = Z(f)$ be a cubic surface, $f_S = t^3 - f$ the cubic form defining its associated cyclic cubic $3$-fold $T = Z(f_S)$, and $\Omega$ the standard volume form on $\PP^4$. The map 
    \begin{align*}
        \CC[z_0 ,z_1,z_2,z_3, t]_{(1)} &\longrightarrow H^{2,1}(T)\\
        P &\longmapsto \mathrm{Res}_T \left(\frac{P \Omega}{f_S^2}\right)
    \end{align*}
    is an isomorphism of vector spaces, under which the line $\CC\langle t \rangle$ maps to $H_{\bar{\omega}}^{2,1}(T)$. 
\end{proposition}

If we additionally assume that $f$ defined an symmetric cubic surface, the cyclic cubic $3$-fold $T$ also admits an $S_4$ symmetry. Thus the cohomology of $T$ inherits the structure of an $S_4$-representation, which we seek to determine.

\begin{lemma}\label{ThreefoldRepresentationLemma}
For any symmetric cyclic cubic 3-fold $T$, the cohomology group $H^{2,1}(T)$ is isomorphic as an $S_4$-representation to $\CC \oplus V$, where $V$ is the standard $S_4$-permutation representation on $\CC^4$.
\end{lemma}
\begin{proof}
The meromorphic differential forms 
$$\left\langle \frac{z_0\Omega}{f_S^2},\frac{z_1\Omega}{f_S^2},\frac{z_2\Omega}{f_S^2},\frac{z_3\Omega}{f_S^2},\frac{t\Omega}{f_S^2}\right\rangle$$
give us a basis for $H^{2,1}(T)$. From this we can explicitly compute the induced structure on $H^{2,1}(T)$ as an $S_4$-representation. Recall that $S_4$ acts by linear permutation automorphisms on $z_0,\dots, z_3$ and acts trivially $t$. Moreover, $\sigma^* \Omega = \Omega$ and $\sigma^* f_S = f_S$ for all $\sigma \in S_4$ since the symmetric group leaves invariant the cubic form $f_S$ and the volume form $\Omega$. The claim follows. 
\end{proof}

\subsection{Equivariant framings and the local period map}
For every cyclic cubic threefold $T$, we have that 
$$H_{\bar{\omega}}^3(T) \cap H^{1,2}(T) = H_{\bar{\omega}}^{1,2}(T)$$ 
is a positive hyperplane in the signature $(4,1)$-space $H_{\bar{\omega}}^3(T)$. Since $S_4$ acts on $T$ by holomorphic automorphisms and commutes with the deck group $\langle \tau \rangle$, the induced action on $H_{\bar{\omega}}^3(T)$ must act trivially on the line $H_{\bar{\omega}}^{2,1}(T)$ defining the period data. After complex conjugating, \autoref{ThreefoldRepresentationLemma} tells us that $H_{\bar{\omega}}^{1,2}(T)$ is isomorphic, as an $S_4$-representation, to the standard permutation representation. The hyperplane $H_{\bar{\omega}}^{1,2}(T)$ is then uniquely determined as an $S_4$-representation by the 1-dimensional trivial $S_4$-representation $\CC \subset H_{\bar{\omega}}^{1,2}(T)$. 

The ambient signature $(4,1)$-space $\CC^{4,1} \cong H_{\bar{\omega}}^3(T)$ is where periods of cubic surfaces $S$ live in. To refine our period data equivariantly, we will use the fixed locus $H_\omega^3(T)_{1}$ of the $S_4$-action on $H_\omega^3(T)$ to define the period domain of symmetric cubic surfaces. By the above discussion, $H_\omega^3(T)_{1}$ is a signature $(1,1)$ complex inner product space.

Let $T$ be a $S_4$-invariant cyclic cubic threefold associated to the symmetric cubic surface $S$, and let $\sigma_T : S_4 \times H^3(T,\ZZ) \to H^3(T,\ZZ)$ be the $S_4$-action induced on the $\ZZ[\omega]$-module $H^3(T,\ZZ)$. Set $\Lambda$ equal to the unique $\ZZ[\omega]$-lattice of signature $(4,1)$ abstractly isomorphic to $H^3(T,\ZZ)$, and $\sigma$ an $S_4$-action on $\Lambda$ abstractly isomorphic to the action $\sigma_T$ on $H^3(T,\ZZ)$. 

\begin{definition}
    An \textit{equivariant framing} is a pair $(S,\lambda)$ of a symmetric cubic surface $S$ and a framing
    $$\lambda : (H^3(T,\ZZ), \sigma_T) \xrightarrow[]{\sim} (\Lambda, \sigma)$$
    which sends the action $\sigma_T$ on $H^3(T,\ZZ)$ to the action $\sigma$ on $\Lambda$. Two equivariant framings $(S_1, \lambda_1)$ and $(S_2, \lambda_2)$ are equivalent if there exists some $g \in N_{\PGL(4,\CC)}(S_4)$ such that $\lambda_1 = \tilde{g}^* \lambda_2$. Let $\mathcal{G}$ denote the space of equivalence classes of equivariantly framed symmetric cubic surfaces $(S,\lambda)$.
\end{definition}

\begin{proposition}[{\cite[Proposition 4.2]{YuZheng}}]
    The space $\mathcal{G}$ is a complex manifold.
\end{proposition}
Let $\Lambda_{\CC,1} \subset \Lambda_\CC = \Lambda \otimes_{\ZZ[\omega]} \CC$ denote the fixed locus of the $S_4$-action $\sigma$ on $\Lambda_\CC$.

\begin{definition}
    The \textit{symmetric period domain} $\mathbb{D}$ associated to the moduli of equivariantly marked symmetric cubic surfaces $\mathcal{G}$ is the Hermitian symmetric domain
    $$\mathbb{D} = \PP\{x \in \Lambda_{\CC,1} \cong \CC^{1,1} : h(x,\overline{x}) <0\}.$$
    Clearly $\mathbb{D} \cong \CC\HH^1$, the complex hyperbolic line. Equivalently, $\mathbb{D}$ is the real hyperbolic plane.
\end{definition}

The following diagram contains most of the spaces of interest. The main content of this section will be showing injectivity of top left horizontal map $\mathcal{G} \to \mathcal{F}$. The right column of horizontal period maps are injective by \cite{AllcockCarlsonToledo}. The remaining horizontal maps in the left column are injective by definition and the pullback construction.
\begin{center}
    \begin{tikzcd}
\text{framed moduli} & {\mathcal{G}} \arrow[d] \arrow[r] & {\mathcal{F}} \arrow[d] \arrow[r] & \CC\HH^4 \arrow[d]                        \\
\text{marked moduli} & {\widehat{\mathcal{S}}} \arrow[d] \arrow[r] & {\widehat{\mathcal{M}}} \arrow[d] \arrow[r] & \hat{\Gamma}'\backslash\CC\HH^4 \arrow[d] \\
\text{moduli}        & {\mathcal{S}} \arrow[r]           & {\mathcal{M}} \arrow[r]           & \hat{\Gamma}\backslash\CC\HH^4           
\end{tikzcd}
\end{center}

\begin{proposition}\label{LocalPeriodMapInjective}
    Let $\mathcal{H}^{S_4}  =\mathbb{D} \cap \mathcal{H}$ denote the symmetric discriminant locus in the period domain. The natural map $\mathcal{G} \to \mathcal{F}$ is injective. Thus the local symmetric framed period mapping $\widetilde{\mathcal{P}}: \mathcal{G} \to \mathbb{D}$ injective, and moreover is an open embedding onto its image $\mathbb{D} - \mathcal{H}^{S_4}$. Moreover, its extension to the stable locus $\mathcal{G}^{s}$ is surjective.
\end{proposition}
\begin{proof}
    Suppose we are given two equivariantly framed symmetric cubic surfaces $(S_1 ,\lambda_1),(S_2,\lambda_2) \in \mathcal{G}$ that map to the same point in $\mathcal{F}$. Then there exists an linear isomorphism of varieties $g : S_1 \xrightarrow[]{\sim} S_2$ along with a unique up to deck transformations lift $\tilde{g} : T_1 \xrightarrow[]{\sim} T_2$ satisfying
    $$\tilde{g}^* = \lambda_1^{-1} \circ \lambda_2 : H^3(T_2 , \ZZ) \to H^3(T_1,\ZZ),$$
    which is an isometry of $\ZZ[\omega]$-lattices. Since $\lambda_1$ and $\lambda_2$ are compatible with the $S_4$-action, so is $\tilde{g}^*$. By \cite[Theorem 1.1]{ZhengOrbifold}, the equivariantly framed cubic surfaces $(S_1,\lambda_1)$ and $(S_2,\lambda_2)$ represent the same point in $\mathcal{G}$. This proves injectivity of the map $\mathcal{G} \to \mathcal{F}$. Commutativity of the diagram
    \begin{center}
        \begin{tikzcd}
{\mathcal{G}} \arrow[d] \arrow[r, hook] & {\mathcal{F}} \arrow[d, hook] \\
\mathbb{D} \arrow[r, hook]                    & \CC\HH^4                     
\end{tikzcd}
    \end{center}
    then implies that the local symmetric framed period map $\mathcal{G} \to \mathbb{D}$ is injective. Moreover, since the derivative of the period map $\widetilde{\mathcal{P}} : \mathcal{F} \to \CC\HH^4$ is injective everywhere, so is the derivative of $\widetilde{\mathcal{P}}: \mathcal{G} \to \mathbb{D}$. Thus the local symmetric framed period map induces a diffeomorphism onto its image $\mathbb{D} - (\mathbb{D} \cap \mathcal{H}) = \mathbb{D} - \mathcal{H}^{S_4}$. Since the local period map on the stable framed moduli space $\mathcal{F}^s$ has image $\CC\HH^4$ \cite[Theorem 3.17]{AllcockCarlsonToledo}, the image of the local symmetric period map $\mathcal{G}^s \to \DD$ is surjective, thereby proving the claim.
\end{proof}

\subsection{The global period map}

Since $S_4$ acts on any symmetric cubic surface and its associated cyclic cubic $3$-fold $T$, it embeds into the arithmetic group $\mathrm{U}(4,1,\ZZ[\omega])$ via its action on $H^3(T,\ZZ)$ (note that this map is injective, since $S_4$ embeds into the mod $3$ reduction $\mathrm{PO}(4,1,\FF_3)$ via its action on the $27$ lines). To study how the local period map descends to yield a uniformization of the symmetric moduli space by the global period map, we must determine the normalizer of $S_4$ in a few groups of interest.

\begin{proposition}\label{S4normU41}
    There is an isomorphism of groups $N_{\mathrm{U}(4,1)}(S_4) \cong \mathrm{U}(1,1) \times (\mathrm{U}(1) \cdot S_4)$.
\end{proposition}
\begin{proof}
    Recall that $S_4$ acts on $\CC^{4,1}$ by the standard permutation representation on the positive $4$-space and the negative $1$-space. This splits the space into a direct sum of irreducible representations $W \oplus \CC \oplus \CC$, where $W$ denotes the irreducible $S_4$-representation of dimension $3$, and the sum of the two trivial representations form a signature $(1,1)$-space. By \autoref{prop:S4normalizercentral}, it suffices to determine what the centralizer of $S_4$ is within $\mathrm{U}(4,1)$. By centrality, we can deal with the $3$-dimensional factor and $(1,1)$-factor individually. 
    
    Schur's lemma tells us that the $S_4$-centralizer acts by scalars on $W$, and thus is isomorphic to a copy of $\mathrm{U}(1)$ acting on $W$. On the signature $(1,1)$-factor, the $S_4$-action is trivial, and thus every element of $\mathrm{U}(1,1)$ arises at an automorphism of the representation $\CC\oplus \CC$. Thus the normalizer is the product of the normalizers on each factor, proving that $N_{\mathrm{U}(4,1)}(S_4) \cong \mathrm{U}(1,1) \times (\mathrm{U}(1) \cdot S_4)$.
\end{proof}

\begin{proposition}
    There is an isomorphism of groups $N_{\mathrm{U}(1)\times \mathrm{U}(4)} (S_4) \cong \mathrm{U}(1)\times (\mathrm{U}(1) \cdot S_4)$, where the second $\mathrm{U}(1)$ acts on the permutation representation $V$ by scalars.
\end{proposition}
\begin{proof}
    Using \autoref{prop:S4normalizercentral}, we need only determine the centralizer of $S_4$ to generate the normalizer. Yet again, these are the unitary scalar matrices.
\end{proof}

\begin{proposition}\label{S4normArithmetic}
    The group $\Gamma = N_{\mathrm{U}(4,1,\ZZ[\omega])}(S_4)$ is naturally an arithmetic subgroup of $N_{\mathrm{U}(4,1)}(S_4)$. Moreover, we have an isomorphism $\Gamma \cong \Aut(\mathrm{diag}(4,-1),\ZZ[\omega]) \times (\langle -\omega \rangle \cdot S_4)$ 
\end{proposition}
\begin{proof}
    For arithmeticity, see {\cite[Appendix A]{YuZheng}}. As before, \autoref{prop:S4normalizercentral} tells us that it suffices to determine the centralizer in $\mathrm{U}(4,1,\ZZ[\omega])$, which we shall do on each factor of the $S_4$-representation. The Eisenstein lattice  $\ZZ[\omega]^{5} \subset \CC^{4,1}$ intersects the $S_4$-representations $W$ and $\CC\oplus \CC$ in rank $2$ and $3$ Eisenstein lattices, respectively. The centralizer is then a subgroup of the product of the automorphism group of these lattices. These symmetries must be automorphisms of $\CC^{4,1}$ which preserve the whole Eisenstein lattice. 
    
    By \autoref{S4normU41}, the $S_4$-centralizer of the rank $3$ sublattice must be scalars $\theta \in \mathrm{U}(1)$ which preserves $\ZZ[\omega]$; this subgroup is $\langle - \omega \rangle$. A standard calculation tells us that generators for this rank $2$ sublattice are given by $(1,1,1,1,0)$ and $(0,0,0,0,1)$, thus the signature $(4,1)$ form restricts to the bilinear form $\mathrm{diag}(4,-1)$. Since $S_4$ acts trivially on this rank $2$ sublattice, the full automorphism group $\Aut(\mathrm{diag}(4,-1),\ZZ[\omega])$ which preserves the Eisenstein lattice constitutes the centralizer on this factor. This proves the claim.
\end{proof}

From these calculations, we can conclude the following which is an application of a well-known fact about locally symmetric varieties \cite[Proposition A.1]{YuZheng}:

\begin{proposition}\label{TotGeodesicEmbedding}
    Let $\hat{G} = \mathrm{U}(4,1)$, $\hat{K} = \mathrm{U}(4)\times \mathrm{U}(1)$, and $\hat{\Gamma} = \mathrm{U}(4,1,\ZZ[\omega])$. Set $G = N_{\hat{G}}(S_4)$, $K = N_{\hat{K}}(S_4)$, and $\Gamma = N_{\hat{\Gamma}}(S_4)$. The holomorphic totally geodesic embedding $G/K \to \hat{G}/ \hat{K}$ descends to a generically injective finite normalization $$\Gamma \backslash G / K \to \hat{\Gamma} \backslash \hat{G} / \hat{K}.$$
\end{proposition}

Following \cite[Proposition 4.10]{YuZheng}, we can prove this more precise version of \autoref{IntroMainThm1}:

\begin{theorem}\label{MainTheoremActually}
    The local period map $\widetilde{\mathcal{P}} : \mathcal{G} \to \mathbb{D}$ descends to isomorphisms $\mathcal{S} \cong P\Gamma \backslash (\DD-\mathcal{H}^{S_4})$ and $\mathcal{S}^s \cong P\Gamma \backslash \mathbb{D}$. Moreover, this is an isomorphism of analytic orbifolds compatible with the uniformization of the moduli of cubic surfaces, so the totally geodesic embedding
    $$P\Gamma \backslash \DD \to P\hat{\Gamma} \backslash \CC\HH^4$$
    is a modular embedding of locally symmetric orbifolds. This map compatibly extends to an isomorphism of the semistable symmetric moduli space $\mathcal{S}^{ss} \cong \overline{P\Gamma \backslash \mathbb{D}}$, where the latter space denotes the Satake compactification of the arithmetic quotient.
\end{theorem}
\begin{proof}
    We will first show that the map $\widetilde{\mathcal{P}}$ descends to a well-defined map $\mathcal{P}([S]) = [\PP(H_{\bar{\omega}}^{2,1}(T))]$. Let $f_1,f_2 \in \mathcal{V}^{\sm}$ be two smooth symmetric cubic forms with equivariant framings $\lambda_1,\lambda_2$ of their associated cubic $3$-folds $T_1,T_2$. Suppose there is some $g \in N_{\PGL(4,\CC)}(S_4)$ such that $g(f_1) = f_2$. This induces the $\ZZ[\omega]$-isometry
    $$\tilde{g}^* : H^3 (T_2,\ZZ) \to H^3(T_1,\ZZ).$$
    We will show that $\gamma = \lambda_1\circ\tilde{g}^*\circ\lambda_2^{-1} \in \Gamma$. Since $g \in N_{\PGL(4,\CC)}(S_4)$, $g \sigma g^{-1}  = \sigma' \in S_4$, thus we have a commutative diagram  
    \begin{center}
        \begin{tikzcd}
\Lambda \arrow[d, "\sigma'"'] \arrow[r, "\lambda_2^{-1}"] & {H^3(T_2,\ZZ)} \arrow[d, "{\sigma' }^*"] \arrow[r, "\tilde{g}^*"] & {H^3(T_1,\ZZ)} \arrow[d, "\sigma*"] \arrow[r, "\lambda_1"] & \Lambda \arrow[d, "\sigma"] \\
\Lambda \arrow[r, "\lambda_2^{-1}"']                      & {H^3(T_2,\ZZ)} \arrow[r, "\tilde{g}^*"']                          & {H^3(T_1,\ZZ)} \arrow[r, "\lambda_1"']                     & \Lambda                    
\end{tikzcd}
    \end{center}
    Thus, as automorphisms of $\Lambda$, $\sigma' = \gamma^{-1} \sigma \gamma$, proving that $\gamma \in \Gamma$. This proves the map $\mathcal{P}$ is a well-defined and yields a commutative diagram
    \begin{center}
        \begin{tikzcd}
\mathcal{G} \arrow[d] \arrow[r, "\widetilde{\mathcal{P}}"] & \DD \arrow[d]         \\
\mathcal{S} \arrow[r, "\mathcal{P}"']                      & P\Gamma\backslash \DD
\end{tikzcd}
    \end{center}
    The Riemann extension theorem tells us that $\mathcal{P}$ extends uniquely to the stable locus $\mathcal{S}^s$. By commutativity of this diagram and \autoref{LocalPeriodMapInjective}, the global period map $\mathcal{P}: \mathcal{S} \to P\Gamma \backslash (\DD-\mathcal{H}^{S_4})$ and its stable extension $\mathcal{S}^s \to P\Gamma \backslash \DD$ are surjective. We shall now show that $\mathcal{P}$ is injective.

    Let $(S_1,\lambda_1),(S_2,\lambda_2) \in \mathcal{G}$ be two equivariantly framed symmetric cubic surfaces with associated cubic forms $f_1,f_2$, and suppose their periods represent the same point in $P\Gamma\backslash\DD$. Then there exists some $\gamma \in \Gamma$ such that $\gamma \cdot \lambda_1(H_{\bar{\omega}}^{2,1} (T_1)) =\lambda_2 (H_{\bar{\omega}}^{2,1} (T_2))$. Thus the map $\lambda_2^{-1} \circ \gamma \circ \lambda_1 : H^3(T_1,\ZZ) \to H^3 (T_2,\ZZ)$ preserved preserves the Eisenstein Hodge structures. By \cite[Theorem 1.1]{ZhengOrbifold}, there exists some $g \in \PGL(4,\CC)$ such that $g(f_2) = f_1$ and $\tilde{g}^* = \lambda_2^{-1} \circ \gamma \circ \lambda_1$. To prove injectivity of $\mathcal{P}$, we want to show that $g \in N_{\PGL(4,\CC)}(S_4)$. 
    
    For any $\sigma \in S_4$ acting on $S_1 = Z(f_1)$, we have that $g^{-1} \sigma g$ acts on $S_2 = Z(f_2)$, which induces the following on the cohomology of the associated cyclic cubic $3$-folds:
    $$(\tilde{g}^{-1} \sigma \tilde{g})^* = \tilde{g}^* \sigma^* (\tilde{g}^{-1})^* = (\lambda_2^{-1} \gamma \lambda_1) (\lambda_1^{-1} \sigma^* \lambda_1) (\lambda_1^{-1} \gamma^{-1} \lambda_2) = \lambda_2^{-1} \gamma \sigma^* \gamma^{-1} \lambda_2.$$
    Since $\gamma \in \Gamma$, we have that $\gamma \sigma^* \gamma^{-1} \in S_4$ as an automorphism of cohomology. Again by \cite[Theorem 1.1]{ZhengOrbifold}, we have that $g \sigma g^{-1} \in S_4$, proving that $g \in N_{\PGL(4,\CC)} (S_4)$. Modular compatibility of the totally geodesic embedding is a consequence of \autoref{TotGeodesicEmbedding}. This proves the global period map satisfies the claimed properties.

    By \cite[Theorem 8.2]{AllcockCarlsonToledo}, the period map extends to the semistable locus for the total moduli space $\mathcal{M}^{ss} \to \overline{P\hat{\Gamma} \backslash\CC\HH^4}$ and sends the unique semistable non-stable point to the unique boundary point of the Satake compactification. Since the embedding of locally symmetric orbifolds $P\Gamma \backslash \DD \to P\hat{\Gamma} \backslash \CC\HH^4$ is modular, the extension to their Satake compactifications is modular, and thus the tricuspidal point on $\mathcal{S}^{ss}$ maps to the unique boundary point of $\overline{P\Gamma \backslash \DD}$, as claimed.
\end{proof}

Now that we have successfully uniformized the moduli space of symmetric cubic surfaces, we will begin our study of the monodromy group associated to the cover $\widetilde{\mathcal{S}} \to \mathcal{S}$, where $\widetilde{S}$ is the moduli of symmetric cubic surfaces equipped with a line. The following calculation implies \autoref{IntroMainThm2}:

\begin{proposition}\label{prop:HodgeMonodromyRestriction}
    Consider the group homomorphism $\mathrm{U}(4,1,\ZZ[\omega]) \to \mathrm{PO}(4,1,\FF_3)$ induced by the map $\omega \mapsto 1$. The subgroup corresponding to $\Gamma = N_{\mathrm{U}(4,1,\ZZ[\omega])} (S_4)$ has image $K_4 \times S_4$ in $\mathrm{PO}(4,1,\FF_3)$.
\end{proposition}
\begin{proof}
    We appeal to the isomorphism explicitly traced out in the proof of \autoref{S4normArithmetic}, and determine the mod $3$ reduction on each factor. As previously discussed, the $S_4$ factor survives the quotient by its action on the $27$ lines. Since $\omega \mapsto 1$, the $\langle - \omega \rangle$ factor which acted on the rank $3$ lattice by scaling is sent to $\langle -1 \rangle \cong C_2$. Finally, since the quadratic form $\mathrm{diag}(4,1)$ reduces mod $3$ to the quadratic form $\mathrm{diag}(1,-1)$, one can calculate that the group $\Aut(\mathrm{diag}(4,-1),\ZZ[\omega])$ has image isomorphic to $\mathrm{PO}(1,1,\FF_3) \cong C_2$. Thus we've shown the image of $\Gamma$ is isomorphic to $K_4 \times S_4$.
\end{proof}

\begin{proof}[Proof of \autoref{IntroMainThm2}]
    Since the inclusion $\mathcal{S} \to \mathcal{M}$ induces an injection on orbifold fundamental groups $P\Gamma \to P\hat{\Gamma}$, it suffices to determine the image of $\Gamma$ in $\PO(4,1,\FF_3)$. This is carried out in \autoref{prop:HodgeMonodromyRestriction}, and is isomorphic to $K_4\times S_4$.
\end{proof}

One may be tempted to conclude that the monodromy group of the cover of parameter spaces $\widetilde{\mathcal{Y}^{\sm}} \to \mathcal{Y}^{\sm}$ is $K_4 \times S_4$. Indeed, \cite[Section 8]{RealACT} outlines why the monodromy groups associated to the connected components of moduli of real projective cubic surfaces are the image of their fundamental groups in $\PO(4,1,\FF_3)$. However, all that \autoref{prop:HodgeMonodromyRestriction} guarantees is that the monodromy group is contained in $K_4 \times S_4$. Remarkably, this fails to pin down our desired Galois group from purely Hodge-theoretic considerations --- it will be a proper subgroup of $K_4 \times S_4$! Further analysis using equivariant line geometry on cubic surfaces is required, which we carry out in the next section. 

\section{Calculating the monodromy group}\label{sec:monodromy-group}

In this section, we will determine the monodromy group of the cover $\widetilde{\mathcal{Y}^{\sm}} \to \mathcal{Y}^{\sm}$ by a combination of classical, moduli-theoretic, and computational techniques. We begin with the following basic fact:

\begin{proposition} The automorphism group of a cubic surface acts faithfully on its lines.
\end{proposition}
\begin{proof}
    Any automorphism $\varphi$ of a cubic surface $S$ preserves the canonical class $K_S$, and thus extends to the anticanonical embedding of $S$ into $\PP^3$, so $\varphi\in\PGL(4,\CC)$. Thus $\varphi$ sends lines to lines on $S$, and any such $\varphi$ which fixes all $27$ lines must be the identity. 
\end{proof}

\begin{example} For symmetric cubic surfaces, this implies that $S_4 \subseteq W(E_6)$. A priori for different symmetric cubic surfaces we might obtain different conjugacy classes of $S_4$ in $W(E_6)$, however the connectivity of the moduli space of symmetric cubic surfaces guarantees this cannot occur. Thus when we discuss $S_4$ as a subgroup of $W(E_6)$ we are implicitly referring to this specific conjugacy class of subgroups.
\end{example}

\begin{proposition}\label{prop:monodromy-in-NWS4}
The symmetric monodromy group is a subgroup of $N_{W(E_6)}(S_4) \cong K_4 \times S_4$.
\end{proposition}
\begin{proof}
This is immediate by translating \autoref{prop:HodgeMonodromyRestriction} along the exceptional isomorphism $W(E_6)\cong \PO(4,1,\FF_3)$. It can also be proved by leveraging Luna's \'etale slice theorem (c.f. \cite{Luna,PopovVinberg}) to argue there exists a universal deformation space for $S_4$-symmetric cubic surfaces, and therefore via descent, monodromy in the symmetric locus preserves the fiberwise $S_4$-action on a universal family of symmetric cubic surfaces, thereby normalizing $S_4$ in the full monodromy group $W(E_6)$.

The splitting of the short exact sequence
\begin{equation}\label{eqn:SES-aut-monodromy}
\begin{aligned}
    0 \to S_4 \to N_{W(E_6)}(S_4) \to K_4 \to 0.
\end{aligned}
\end{equation}
claimed in the proposition is a computer verifiable computation.
\end{proof}

\subsection{Certified tracking}

In order to gain some insight into the structure of the symmetric monodromy group, we wish to witness the existence of certain elements by lifting explicit loops in the parameter space. In conversations with T. Brysiewicz, working with his Pandora software \cite{Pandora}, we were able to generate strong computational evidence towards the structure of the monodromy group.

Algorithms used in this and related software fall under the umbrella of \textit{homotopy continuation}. This is a key technique in numerical algebraic geometry which deforms a system of polynomial equations along a one-parameter path. One of the primary applications of this technology is conducting explicit monodromy computations.

While homotopy continuation software can generate strong evidence towards a computation, more refined algorithms are needed to turn these computations into proof. At each stage of tracking solutions along a one-parameter path, a guarantee is needed that paths don't collide, and therefore that the computed permutation is indeed correct. These more sophisticated (and time-costly) methods are called \textit{certified tracking algorithms}. Recent work of T. Duff and K. Lee provides algorithms which, among other things, are applicable for certifying computations in monodromy, bridging the gap between computation and proof \cite[Theorem~1]{DuffLee}. In conversations with Lee, their software is able to mathematically certify \autoref{thm:IntroMainThm3}:

\begin{theorem}[Numerical certification]\label{lem:certification}
The monodromy group of lines on symmetric cubic surfaces is isomorphic to the Klein 4-group, and centralizes $S_4$ in $W(E_6)$.
\end{theorem}
\begin{proof}[Setup] To determine the monodromy group, it suffices to obtain generators for $\pi_1(\mathcal{Y}^\sm)$ and then run certified tracking algorithms on their fibers. The discriminant locus of cubic surfaces, when intersected with the symmetric locus, factors into a product of four irreducible components, three of which are real lines and one of which is a cubic polynomial. Taking loops around each of these and lifting gives the desired result. Finally to identify $W(E_6)$ as a subgroup of $S_{27}$, we must label the numerical solutions for lines at our basepoint with explicit equations for the lines on the Fermat cubic surface.
\end{proof}

\begin{remark}[On this monodromy group] \,
\begin{enumerate}
    \item There are further geometric constraints on monodromy that fail to completely pin down the group --- for instance, the tritangent on $S_4$-symmetric cubic surfaces which is stabilized by the dihedral group $D_8$ \cite[Theorem 1.2]{EEG} is fixed pointwise under any path in the smooth symmetric locus. This implies the monodromy group is contained in its stabilizer in $W(E_6)$, a subgroup of order $192$ (see \autoref{prop:monodromy-contained-in-tritangent-stabilizer} below). The intersection of this subgroup with the normalizer of $S_4$ has order $16$, and contains the honest monodromy group properly. However, knowing the stabilizer of this tritangent will be useful when determining the connected components of the symmetric $27$ lines cover $\widetilde{\mathcal{Y}} \to \mathcal{Y}$, which we do in the next section.
    
    \item The short exact sequence \autoref{eqn:SES-aut-monodromy} admits an explanation using the language of stacks in recent work of A. Landi \cite{Alberto}. In this work, Landi resolves various equivariant monodromy problems by defining new stacks of equivariant objects. In particular, he recovers this monodromy calculation with completely orthogonal methods.

    \item After the first draft of this paper appeared, E. Pichon-Pharabod and S. Telen were inspired to approach these ideas with different numerical techniques. In particular they numerically compute the monodromy groups of lines on cubic surfaces with other automorphism groups by certifying their action on cohomology. They recover our \Cref{lem:certification} as part of \cite[Theorem~2]{pichonpharabod2025galoisgroupssymmetriccubic}.
    
\end{enumerate}
\end{remark}

\subsection{The incidence variety of 27 lines over the symmetric locus}

Now that we have determined the monodromy group of the cover $\widetilde{\mathcal{Y}^{\sm}} \to \mathcal{Y}^{\sm}$ is $K_4$, we would like to say more about the topology of the space of symmetric cubic surfaces with a line, namely how this restricted cover splits into connected components. To do this, we will make explicit the action of $K_4$ on $\widetilde{\mathcal{Y}^{\sm}}$. 

We first recall the following result of the first named author, which proves that the $S_4$-action on the $27$ lines is independent of the choice of smooth symmetric cubic surface: 

\begin{theorem}\label{thm:eeg} \cite[Theorem~1.2]{EEG} On any smooth symmetric cubic surface, the 27 lines have orbits
\begin{align*}
    \left[ S_4/C_2^o \right] + \left[ S_4/C_2^e \right] + \left[ S_4/D_8 \right],
\end{align*}
where $C_2^o = (1 \ 2)$ is an odd copy of the cyclic group of order two, and $C_2^e = (1\ 2)(3\ 4)$ is an even copy of the cyclic group of order two.
\end{theorem}

\begin{example} The \textit{Fermat cubic surface} is defined by the symmetric homogeneous form $m_3$. Its 27 lines, with explicit labels and parametric equations, are given in the appendix of this paper (\autoref{data:Fermat-lines}). The lines $\ell_1, \ldots, \ell_{12}$ lie in the $S_4/C_2^o$ orbit, the lines $\ell_{13}, \ldots, \ell_{24}$ lie in the $S_4/C_2^e$ orbit, and the lines $\ell_{25},\ell_{26},\ell_{27}$ form a tritangent which is the $S_4/D_8$ orbit. We refer to these three lines as the \textit{$D_8$-tritangent}.
\end{example}

\begin{remark} Once labels are fixed on the 27 lines, we can construct $W(E_6)$ as a permutation group, given as the adjacency-preserving permutations of the 27 lines. As a subgroup of $S_{27}$ with the labeling of the lines coming from the Fermat cubic surface, the generators for $W = W(E_6)$ are listed in \autoref{data:fermat-W}.
\end{remark}

\begin{proposition} The three lines $\{\ell_{25},\ell_{26},\ell_{27}\}$ lie on every symmetric cubic surface, forming a $D_8$-tritangent. Moreover they are fixed under symmetric monodromy.
\end{proposition}
\begin{proof} Since each symmetric cubic surface is a linear combination of elementary homogeneous symmetric polynomials, it suffices to verify each of these vanishes on the lines in the $D_8$-tritangent, which is a routine computation. 

Since symmetric monodromy is $S_4$-equivariant, the tritangent plane spanned by the lines $\ell_{25},\ell_{26},\ell_{27}$ must be stabilized. Moreover, there is no fourth distinct line incident to any symmetric cubic surface which lies in the $D_8$-tritangent, as this would violate B\'{e}zout's theorem. Any nontrivial deformation of $\ell_{25},\ell_{26},\ell_{27}$ arising from monodromy would yield such a line, and so the lines $\ell_{25},\ell_{26},\ell_{27}$ must be fixed by monodromy within the symmetric locus.
\end{proof}

Observe what this means --- given any loop in the symmetric locus, viewed as an element of $W(E_6) \le S_{27}$, it fixes each of the points 25, 26, and 27. Since $W(E_6)$ acts transitively on ordered tritangents, we can ask what the pointwise stabilizer of a tritangent is in $W(E_6)$, and this will contain our monodromy group.

\begin{proposition}\label{prop:monodromy-contained-in-tritangent-stabilizer}
The symmetric monodromy group is contained in the pointwise stabilizer of a tritangent in $W(E_6)$. This is a group of order 192.
\end{proposition}

By combining our constraints for the symmetric monodromy group arising from uniformization (\autoref{prop:monodromy-in-NWS4}) and from equivariant enumerative geometry (\autoref{prop:monodromy-contained-in-tritangent-stabilizer}), we obtain the following reduction.

\begin{proposition}\label{StabIntersectNormalizer}
The monodromy group is contained in the group of order 16:
\begin{align*}
    \bigcap_{i=25}^{27} \Stab_{W(E_6)}(\ell_i) \cap N_{W(E_6)}(S_4) \cong K_4 \times K_4.
\end{align*}
We give names to these generators. The former is $K_4 = \left\langle \sigma_1,\sigma_2 \right\rangle$, and it is a subgroup of $S_4$. The latter is $K_4 = \left\langle \tau_1,\tau_2 \right\rangle$ and it is not contained in $S_4$. As explicit elements in $W(E_6) \le S_{27}$ they are listed in \autoref{data:gens-K4-x-K4}.
\end{proposition}

\begin{corollary}\label{cor:incidence-var-structure} The incidence variety of 27 lines restricted to the symmetric locus $\widetilde{\mathcal{Y}^{\sm}}$ has 12 connected components. Explicitly as a $K_4$-set, the fiber over any symmetric cubic surface is of the form
\begin{align*}
    6\left[ K_4/C_2 \right] + 3\left[ K_4/e \right] + 3\left[ K_4/K_4 \right].
\end{align*}
\end{corollary}
\begin{proof}
    Having restricted the symmetric monodromy group and concluding that it is $K_4 < W(E_6)$, we can then see explicitly how $K_4$ acts on and stablizes the $27$ lines on the Fermat cubic surface. This splits them into the $12$ families claimed (see \autoref{data:gens-K4-x-K4} for the relevant generators and how they act on the $27$ lines).
\end{proof}

To conclude, we have an equivariant parameter space analog of \autoref{WholeBeauvilleMonodromy}:

\begin{theorem}
    The space $\widetilde{\mathcal{Y}^{\sm}}$ is a naturally a (disconnected) complex manifold. Moreover, the equivariant monodromy representation
    $$\pi_1(\mathcal{Y}^{\sm},X) \to \Aut_{S_4}(H^2 (X,\ZZ),\eta_X)\cong S_4 \times K_4$$
    is \textit{not} surjective and has image isomorphic to the Klein $4$-group $K_4$.
\end{theorem}

\section{Symmetry and monodromy via representation theory}\label{sec:repTheory}

In the previous section, we determined the monodromy group $K_4$ within the Weyl group $W(E_6)$ in terms of how it acts on the lines of the Fermat cubic surface; the $S_4$-orbits of the $27$ lines are given in \autoref{data:Fermat-lines}. The goal of this section is to understand these copies of $S_4$ and $K_4$ in $W(E_6)$ from more traditional representation theoretic viewpoint, via reflection group theory and the projective orthogonal perspective. Informally, we will show that the symmetry group and monodromy group are not visible from purely Coxeter-theoretic considerations.

\subsection{The Weyl group as a reflection group}
To present the Weyl group $W(E_6)$ as a reflection group, we first label the nodes of $E_6$ Dynkin diagram with the generating reflections $s_0,\dots,s_5$:
\begin{center}
    \dynkin[edge length=1.5cm,labels*={s_1,s_0,s_2,s_3,s_4,s_5}]E6
\end{center}
This gives rise to a presentation of the Weyl group of $E_6$ as a Coxeter group:
\begin{align*}
    W(E_6) = \left\langle s_0, \ldots, s_5 | (s_i s_j)^{m_{ij}} = 1 \right\rangle, \quad\quad m_{ij} = \begin{cases} 1 & i=j \\ 3 & s_i,\ s_j\text{ share an edge} \\
    2 & \text{otherwise}\end{cases}
\end{align*}

The following result is very classical, we recap it for the reader's convenience.
\begin{proposition} Any choice of six skew lines gives rise to a presentation of $W(E_6)$ in the form above.
\end{proposition}
\begin{proof}
    The choice of six skew lines determine a marking of the homology of a cubic surface $S$, where each line corresponds to the homology classes of orthogonal $(-1)$-exceptional curves $e_1,\ldots,e_6$ on $S$. These in turn give us a basis of long roots for the $E_6$ lattice $v_0 = h - e_1 - e_2 - e_3, v_j = e_j - e_{j+1}$ for $j = 1,\dots,5$. The intersection form $Q$ on the homology $H_2(S,\mathbb{Z})$ satisfies
\begin{align*}
    Q(h,h) &= 1 \\
    Q(e_i,e_j) &= -\delta_{ij} \\
    Q(h,e_i) &= 0.
\end{align*}
From this it is clear that $Q(v_i,v_i) = -2$ for any $0\le i \le 5$. Then the reflections $s_i$ that generate the Weyl group $W(E_6)$ are realized homologically by
$$s_i(x) = x - \frac{2Q( x,v_i )}{Q(v_i,v_i)} v_i = x + Q(x,v_i)v_i;$$
this is the \textit{geometric representation} of $W(E_6)$.
\end{proof}

\begin{example} If we pick the lines $[\ell_1,\ell_3,\ell_{10},\ell_{11},\ell_{16},\ell_{22}]$, a direct computation gives the six generators of $W(E_6)$ as the following permutations in $S_{27}$:
\begin{center}
\begin{tabular}{c|c}
    $s_0$ &  \texttt{(1,8)(3,6)(9,26)(10,25)(13,21)(20,23)}\\
    $s_1$ & \texttt{(1,3)(2,4)(5,7)(6,8)(14,15)(18,19)} \\
    $s_2$ & \texttt{(2,12)(3,10)(5,27)(6,25)(14,17)(19,24)} \\
    $s_3$ & \texttt{(5,8)(6,7)(9,12)(10,11)(17,20)(21,24)} \\
    $s_4$ & \texttt{(5,14)(7,15)(9,13)(11,16)(17,27)(21,26)} \\
    $s_5$ & \texttt{(13,23)(14,19)(15,18)(16,22)(17,24)(20,21)}.
\end{tabular}
\end{center}
\end{example}

It is a classical computation that there are exactly 72 ways to pick six pairwise skew lines on a cubic surface.

\subsection{Double sixes from the Weyl group}

Given six ordered pairwise skew lines, we obtain an associated subgroup $W(A_5) \le W(E_6)$ by suppressing the node $s_0$, and all of these subgroups are conjugate. We note though, that we can permute the ordering of our six lines -- a natural question to ask is whether such a permutation extends to element of the Weyl group, and if such an extension exists, whether it is unique. The answer to both these questions is yes.

\begin{proposition} Given six ordered skew lines, any automorphism $\sigma$ of them extends uniquely to an adjacency-preserving automorphism of all 27 lines, i.e. an element of $W(E_6)$.
\end{proposition}
\begin{proof} Any permutation of the lines permutes the homology classes $e_1, \ldots, e_6$ accordingly, and in particular will fix the canonical class $K_S = 3h - e_1 - \cdots - e_6$. Therefore by definition it extends to an element of $W(E_6)$. Since its action on the $e_i$'s defines its action on $h$ and therefore on a basis of the homology, such an extension is unique.
\end{proof}

Moreover, we understand this subgroup of $W(E_6)$.

\begin{proposition} Fixing six ordered skew lines, the subgroup of $W(E_6)$ obtained by permuting them is exactly equal to the Weyl group $W(A_5)$ obtained from the presentation coming from the choice of lines.
\end{proposition}
\begin{proof} It suffices to show that each of the generators $s_1, \ldots, s_5$ is contained in this symmetric group. This is immediate, since $s_i$ permutes $e_i$ and $e_{i+1}$ and fixes the other $e_j$.
\end{proof}
There is a unique conjugacy class of subgroup $W(A_5) \le W(E_6)$, and $W(E_6)/W(A_5)$ is a transitive set of order 36. There are, however, 72 unordered choices of six skew lines. This gives us a surjection
\begin{align*}
    \left\{ \text{six skew lines} \right\} \to W(E_6)/W(A_5),
\end{align*}
which is 2-to-1. In particular, six skew lines come in pairs, which give rise to the same copy of $W(A_5)$ in $W(E_6)$. These pairs of six skew lines are what are known as \textit{double sixes}.

In particular a computation shows that, as a $W(A_5)$-set, the set of lines $\left\{ 1, \ldots, 27 \right\}$ decompose into two transitive $W(A_5)$-sets of order six, and a single transitive set of order 15. These are the double six, and the remaining lines, respectively.

\begin{remark}\label{rmk:other-S6}
While $W(A_5)$ is isomorphic to $S_6$ as we have seen, it is abuse of terminology to equate them. There are \textit{two} non-conjugate subgroups of $W(E_6)$ which are isomorphic to $S_6$, the first being our $W(A_5)$ group, and the latter just being another subgroup of $W(E_6)$ which we denote by $S_6$. The latter group can be distinguished via its action on 27 lines --- it acts transitively on 12 lines and transitively on the other 15.
\end{remark}

\subsection{Our groups are not reflection groups}

We can now argue that both the $S_4$ acting on symmetric cubic surfaces and the symmetric monodromy group are \textit{not reflection subgroups of $W(E_6)$}. This is perhaps obvious to those familiar with e.g. \cite{Manivel}, but we can give an elementary argument now with the machinery we have built.

\begin{proposition} The subgroup $S_4 \le W(E_6)$ is not a reflection subgroup --- that is, it is not isomorphic to $W(A_3)$ for a presentation of $W(E_6)$ arising from any choice of six skew lines.
\end{proposition}
\begin{proof} We prove something stronger, namely that $S_4$ is not subconjugate to $W(A_5)$. Indeed suppose towards a contradiction that it was. As we have seen by \cite{EEG}, the action of $S_4$ on the 27 lines decomposes into three $S_4$-sets, of order 12, 12, and 3. If $S_4 \le W(A_5)$, then this action would be restricted from the action of $W(A_5)$ on the set of 27 lines. However the partition of $\left\{ 1, \ldots, 27 \right\}$ into orbits will only ever \textit{refine} under a restricted group action. In particular since $W(A_5)$ has two orbits of size six it cannot restrict to the prescribed $S_4$-action.
\end{proof}

\begin{remark} The action of the \textit{other} $S_6$ from \autoref{rmk:other-S6} does not have this same restriction, and a computation shows that $S_4$ is indeed subconjugate to $S_6$ in $W(E_6)$.
\end{remark}

\begin{remark}\label{rmk:preferred-double-six} $\ $
\begin{enumerate}
    \item Another interesting note is that while $S_4$ is not subconjugate to $W(A_5)$, we have that $W(A_5)$ is nested in a maximal subgroup isomorphic to $W(A_5) \times C_2 \le W(E_6)$. It \textit{is} true that $S_4$ is subconjugate to this maximal subgroup, and moreover the centralizer of $S_4$ in $W(A_5) \times C_2$ is identical to the symmetric monodromy group!
    \item There is actually a \textit{unique} copy of $W(A_5)$ in $W(E_6)$ for which $S_4$ is a subgroup of its maximal supergroup $W(A_5)\times C_2$. This unique copy corresponds to a \textit{preferred double six for symmetric cubic surfaces}. A direct computation shows that this is the unique double six where six skew lines lie in the same $S_4$-orbit.
\end{enumerate}
\end{remark}

\begin{proposition}
    The symmetric monodromy group $K_4 \leq W(E_6)$ is not a reflection subgroup.
\end{proposition}
\begin{proof}
    Suppose for the sake of contradiction that $K_4$ was a reflection subgroup; it would then take on the form of $W(A_1) \times W(A_1)$. Since each of the generators $s_i$ act on the 27 lines as a product of six disjoint transpositions, there are two nontrivial elements of $W(A_1) \times W(A_1)$ that are the product of six disjoint transpositions. However, computations in GAP (using \autoref{data:Fermat-S4}) tell us that the symmetric monodromy group only has one element that is the product of six disjoint transpositions, a contradiction.
\end{proof}

\subsection{Symmetric monodromy in the projective orthogonal groups}
Since we know how the symmetric monodromy group $K_4 < W(E_6)$ acts on the $27$ lines of the Fermat cubic surface $S$, we can explicitly connect this $K_4$ back to the projective orthogonal group by a lengthy homological calculation. We sketch this correspondence now. 

Recall that the set of six skew lines $\{\ell_{1},\ell_{3},\ell_{10},\ell_{11},\ell_{16},\ell_{22}\}$ determine a marking of the homology of $S$, where each line corresponds to the homology classes of orthogonal $(-1)$-exceptional curves $e_1,e_2,e_3,e_4,e_5,e_6$ on $S$. Using \autoref{data:Fermat-S4}, we can calculate how the monodromy group $K_4$ acts on the exceptional $(-1)$-curves, which in turns explicitly determines how $K_4$ acts on the $E_6$ lattice. Then by passing to the root lattice quotient used in the proof of the exceptional isomorphism outlined in \autoref{ExceptionalIsomorphismWE6}, this $K_4$ projects to the symmetric monodromy group $K_4$ inside of $\PO(4,1,\FF_3)$.

It would be interesting to understand how the symmetric monodromy group arises purely by an analyzing its action on the associated symmetric cyclic cubic $3$-folds. This leads us to the following problem:

\begin{problem}
    Determine the symmetric monodromy group $K_4$ as a subgroup $\PO(4,1,\FF_3)$ directly, that is, without reference to the action on the lines or the exceptional isomorphism with $W(E_6)$.
\end{problem}

As Beauville remarks \cite[pg. 19]{BeauvilleBourbaki}, what makes this difficult is that it is unknown how to produce a marking of a cubic surface from a framing of the corresponding cyclic cubic 3-fold. A resolution to this problem would shed further light on symmetric monodromy can be witnessed by Hodge theory, and therefore could be applied to similar equivariant enrichments of classical enumerative problems.

\section{A formula in radicals for lines on a symmetric cubic surface}\label{sec:formulas}

Here we work out an explicit formula in radicals for the lines on an $S_4$-symmetric cubic surface. Upcoming work between the two authors and A. Landi proves that the monodromy group and Galois group of symmetric enumerative problems agree. The work below does not depend upon that result, however this upcoming work leads us to expect the existence of a formula in exactly two radicals for lines an $S_4$-symmetric cubic surface. Phrased differently, for every smooth $S_4$-symmetric cubic surface defined over $\mathbb{Q}$, we expect that the lines will be defined over a Klein-four Galois extension $\mathbb{Q}(\sqrt{\alpha},\sqrt{\beta})$, where $\alpha$ and $\beta$ are defined in terms of explicit formulas in terms of the coefficients describing the cubic surface. Indeed this is true.

\begin{theorem}\label{thm:formulas}
    Given a generic $S_4$-equivariant cubic surface $a m_3 + b m_{21} + c m_{111}$ defined over $\mathbb{Q}$, its lines are all defined over the Klein four Galois extension $K = \mathbb{Q}(\sqrt{\sqrtone}, \sqrt{\sqrttwo})$, where
    \begin{align*}
        \sqrtone &= -{\left(9  a^{3} + 9  a^{2} b - 9  a b^{2} + 7  b^{3} - 3  a^{2} c - 6  a b c - 3  b^{2} c + 4  a c^{2}\right)} {\left(3  a + b - c\right)}\\
        \sqrttwo&=-{\left(3  a + b - c\right)} {\left(a + 3  b + c\right)}.
    \end{align*}
    Explicit formulas for lines in each orbit are given parametrically (as images of $\P^1$ with coordinates $[s:t]$) as follows:
    \begin{enumerate}
        \item[($S_4/C_2^o)$] The lines in this orbit are all the $S_4$-orbits of
        \begin{align*}
            \left[s + t : -s+t : \frac{9a^2-b^2-(3a-b)c + \sqrt{\sqrtone}}{6ab+2b^2 - 3(a+b)c+c^2}t : \frac{9a^2-b^2-(3a-b)c - \sqrt{\sqrtone}}{6ab+2b^2 - 3(a+b)c+c^2}t\right]
        \end{align*}
        \item[$(S_4/C_2^e)$] The lines in this orbit are all the $S_4$-orbits of
        \begin{align*}
            \left[\frac{a-b-c+ \sqrt{\sqrttwo}}{2(a+b)}s + t : \frac{a-b-c+ \sqrt{\sqrttwo}}{2(a+b)}s - t : s + \frac{X + Y \sqrt{\sqrttwo}}{2\sqrt{\sqrtone}} t : s - \frac{X - Y \sqrt{\sqrttwo}}{2\sqrt{\sqrtone}} t\right],
        \end{align*}
        where
        \begin{align*}
            X &=-9a^2 - 6ab - b^2 + 2(3a + b)c - c^2,\\
            Y &=3a - 3b + c.
        \end{align*}
        \item[$(S_4/D_8)$] The lines in this orbit are all the $S_4$-orbits of
        \[
            [s: -s : t : -t].
        \]
    \end{enumerate}
\end{theorem}

The remainder of this section is devoted to proving this theorem.

\begin{remark}
The discriminant locus for $S_4$-symmetric cubic surfaces is
\small
\begin{align*}
    \Delta &= (a + 3b + c) (3a - 3b + c)^{10} (3a + b - c)^{9} (9a^3 + 9a^2b - 9ab^2 + 7b^3
- 3a^2c - 6abc - 3b^2c + 4ac^2)^4.
\end{align*}
\normalsize
We recognize that the lines above can fail to exist only when the discriminant vanishes.
\end{remark}

The lines in the $D_8$ tritangent are common to all symmetric cubic surfaces, so it suffices to study lines with cyclic isotropy. For the following calculations, fix representative generators for the odd and even conjugacy classes of $C_2$ in $S_4$: we pick $C_2^o = \langle (1\ 2)\rangle$ and $C_2^e = \langle (1 \ 2) (3 \ 4)\rangle$; all other lines with conjugate isotropy can be obtained by taking $S_4$-orbits. We begin with the lines with odd $C_2$ isotropy.

\begin{proposition}
    If $\ell$ is a line with isotropy exactly equal to $(1\ 2)$, then it passes through the point $[1:-1:0:0]$ and exactly one point on the plane $x=y$.
\end{proposition}
\begin{proof}
    Given a line $\ell$ with isotropy group exactly equal to $\langle (1\ 2)\rangle \le S_4 \le \PGL_4$, it is fixed under this action of the cyclic group of order two. Since this action cannot reverse orientation on the line (which is topologically a 2-sphere), it must be a rotation and hence has at least two fixed points. We compute that the fixed locus of $\P^3$ under a single transposition is a plane and a point, namely
    \begin{align*}
        (\P^3)^{C_2^o} = V(x-y) \cup \{[1:-1:0:0]\}.
    \end{align*}
    This is because if a point $[x:y:z:w]$ is fixed pointwise, then we have that
    $(x,y,z,w) = (\lambda y, \lambda x, \lambda z, \lambda w)$ for some $\lambda$. It is clear to see that $\lambda^2=1$, yielding the two possibilities above.

    Therefore to conclude the proposition, it suffices to argue that $\ell$ cannot be contained in the plane $V(x-y)$. If it was, we could act via $(3\ 4)$ and obtain another line $\ell'$ also contained in the plane. Observe $\ell'$ must be distinct from $\ell$ since otherwise this would violate the isotropy assumption. This implies that $V(x-y)$ would be a tritangent plane to our symmetric cubic surface, and there would necessarily exist a third line $\ell''$ on this plane which is also on the cubic surface. Acting via $(3\ 4)$ on these three lines, we see one must be fixed, implying its isotropy contains the non-normal Klein four group $\langle(1\ 2), (3\ 4)\rangle$. That line necessarily has a $D_8$ isotropy subgroup. However none of the lines in the $D_8$ tritangent lie on the plane $V(x-y)$, a contradiction.
\end{proof}

The formula for lines with odd $C_2$ isotropy follow by solving symbolically for the point of intersection for a line $\ell$ on a symmetric cubic surface and the plane $V(x-y)$.

For the lines with even isotropy, we obtain the following:

\begin{proposition}\label{prop:solving-for-even-lines}
Any line whose isotropy is equal to $(1\ 2)(3\ 4)$ intersects a line in the $D_8$-tritangent, and also intersects one of the following two points:
\begin{align*}
    P^+ &= \left[\frac{a-b-c+ \sqrt{\sqrttwo}}{2(a+b)} : \frac{a-b-c+ \sqrt{\sqrttwo}}{2(a+b)} : 1 : 1\right] \\
    P^- &= \left[\frac{a-b-c- \sqrt{\sqrttwo}}{2(a+b)} : \frac{a-b-c- \sqrt{\sqrttwo}}{2(a+b)} : 1 : 1\right].
\end{align*}
\end{proposition}
\begin{proof}
Again, $\ell$ admits at least two fixed points. These points must lie on the union of the two skew lines
\begin{align*}
    (\P^3)^{C_2^e} = V(z_0-z_1,z_2-z_3) \cup V(z_0+z_1,z_2+z_3).
\end{align*}
We see that $\ell$ cannot be equal to the latter line since that line has isotropy $D_8$, and the former line cannot lie on our cubic surface as it intersects the point $[1:1:-1:-1]$, which lies on the intersection of the other two lines in the $D_8$-tritangent. Therefore $\ell$ intersects both lines. To obtain the formula in the theorem, we solve for the point of intersection $[1:-1:\lambda:\lambda]$ between a line passing through $P^+$. The line passing through $P^-$ can be solved similarly, or obtained by permuting coordinates.
\end{proof}

\begin{remark}
    In the proof of \Cref{prop:solving-for-even-lines}, solving symbolically over the function field $\mathbb{Q}(a,b,c)$, it is not directly obvious that the quantity $\lambda$ is an element of $\mathbb{Q}(a,b,c)(\sqrt{\alpha},\sqrt{\beta})$. This is because $\lambda$ contains a term of the form $\sqrt{t_0 + t_1 \sqrt{\beta}}$ for some quantities $t_0,t_1\in \mathbb{Q}(a,b,c)$. We can argue, however, that this specific quantity $t_0 + t_1 \sqrt{\beta}$ is in fact a square in $\mathbb{Q}(a,b,c)(\sqrt{\beta})$ and indeed $\lambda \in \mathbb{Q}(a,b,c)(\sqrt{\alpha},\sqrt{\beta})$. This observation is necessary to produce the nice forms of the lines with even isotropy.
\end{remark}

\begin{remark}[On tritangents and computation] We have noted in \Cref{rmk:preferred-double-six} that the lines in the orbit $S_4/C_2^o$ form a double six configuration. In particular it is easy, once we have formulas for these lines, to find three skew lines in the orbit. Once we have these in hand, one could solve directly for the equations of all the other lines using the methods in \cite{three_skew}. This is an alternative and equivalent way to obtain the lines with even isotropy, however we have found it computationally easier to solve for them directly as in \Cref{prop:solving-for-even-lines}.
\end{remark}

\begin{example}[Explicit lines on a symmetric cubic surface] When $(a,b,c) = (4,-3,1)$, for instance, we obtain all $27$ lines as the orbit of the following three lines:
\begin{align*}
    &\left[s+t : -s + t : \left(-\frac{1}{14}\sqrt{178} - \frac{9}{4}\right) t : \left(-\frac{1}{14}\sqrt{178} + \frac{9}{4}\right) t\right], \\
    &\left[(2\sqrt{2}+3)s + t : (2\sqrt{2} + 3)s - t : s + \frac{11\sqrt{2}-8}{\sqrt{178}}t : s - \frac{11\sqrt{2}-8}{\sqrt{178}}t\right],\\
    &[s:-s:t:-t].
\end{align*}    
\end{example}

We conclude by providing an algebraic reproof of \cite[1.3]{EEG}.

\begin{theorem}
    A real smooth cubic surface can contain only 3 or 27 real lines.
\end{theorem}
\begin{proof}
    It is clear that when $\alpha,\beta>0$, we obtain 27 real lines via the formulas in \Cref{thm:formulas}. When $\alpha<0$ none of the lines in the odd or even orbits are defined over $\mathbb{R}$. In the lines for the even orbit, we can scale through by $2\sqrt{\alpha}$, and see that the first two entries involve a $\sqrt{\alpha\beta}$ term, while the last two entries involve a $\sqrt{\beta}$ term. These are only defined over $\mathbb{R}$ when $\alpha>0$ and $\beta>0$, in which case all 27 lines are defined over $\mathbb{R}$.
\end{proof}

\appendix
\section{Data tables}\label{sec:Appendix}

We record some of the line geometry data associated to the Fermat cubic surface.

\subsection{All about the Fermat}\label{subsec:fermat-data}
\begin{data}\label{data:Fermat-lines} The 27 lines $\ell_i$ on the Fermat can be labeled and grouped according to their $S_4$-orbits as follows:
\begin{center}
    \begin{tabular}{|l | l |}
    \hline
    $i$ & $\ell_i$  \\
    \hline
1 & $[w, -w, z, \zeta\cdot z]$  \\
2 & $[w, -w, z, \zeta^5\cdot z]$  \\
3 & $[w, \zeta\cdot w, z, -z]$  \\
4 & $[w, \zeta^5\cdot w, z, -z]$  \\
5 & $[w, z, \zeta\cdot w, -z]$  \\
6 & $[w, z, \zeta^5\cdot w, -z]$\\
7 & $[w, z, -w, \zeta\cdot z]$  \\
8 & $[w, z, -w, \zeta^5\cdot z]$ \\
9 & $[w, z, -z, \zeta\cdot w]$  \\
10 & $[w, z, -z, \zeta^5\cdot w]$  \\
11 & $[w, z, \zeta\cdot z, -w]$ \\
12 & $[w, z, \zeta^5\cdot z, -w]$  \\
\hline
\end{tabular}
\quad   
\begin{tabular}{|l | l |}
    \hline
    $i$ & $\ell_i$  \\
    \hline
13 & $[w, \zeta\cdot w, z, \zeta\cdot z]$ \\
14 & $[w, \zeta\cdot w, z, \zeta^5\cdot z]$  \\
15 & $[w, \zeta^5\cdot w, z, \zeta\cdot z]$  \\
16 & $[w, \zeta^5\cdot w, z, \zeta^5\cdot z]$  \\
17 & $[w, z, \zeta\cdot w, \zeta\cdot z]$  \\
18 & $[w, z, \zeta\cdot w, \zeta^5\cdot z]$ \\
19 & $[w, z, \zeta^5\cdot w, \zeta\cdot z]$ \\
20 & $[w, z, \zeta^5\cdot w, \zeta^5\cdot z]$  \\
21 & $[w, z, \zeta\cdot z, \zeta\cdot w]$ \\
22 & $[w, z, \zeta^5\cdot z, \zeta\cdot w]$  \\
23 & $[w, z, \zeta\cdot z, \zeta^5\cdot w]$ \\
24 & $[w, z, \zeta^5\cdot z, \zeta^5\cdot w]$  \\
\hline
\end{tabular}
\quad
    \begin{tabular}{|l | l |}
    \hline
    $i$ & $\ell_i$  \\
    \hline
25 & $[w, -w, z, -z]$ \\
26 & $[w, z, -w, -z]$ \\
27 & $[w, z, -z, -w]$ \\
\hline
\end{tabular}
\end{center}
\end{data}

\begin{data}\label{data:fermat-W} Given the labeling of the lines on the Fermat as in \autoref{data:Fermat-lines}, the Galois group $W(E_6)$ can be maniputaled in GAP (\cite{GAP4}) by
\begin{verbatim}
G := SymmetricGroup(27);
W:= Subgroup(G,[
(13,23)(14,19)(15,18)(16,22)(17,24)(20,21),
(5,14)(7,15)(9,13)(11,16)(17,27)(21,26),
(2,6)(4,8)(5,19)(7,18)(9,23)(11,20)(12,25)(16,21)(22,26)(24,27),
(5,8)(6,7)(9,12)(10,11)(17,20)(21,24),
(3,4)(5,10)(6,9)(7,12)(8,11)(13,15)(14,16)(17,24)(18,23)(19,22)(20,21)(26,27),
(1,2)(5,9)(6,10)(7,11)(8,12)(13,14)(15,16)(17,21)(18,22)(19,23)(20,24)(26,27)
]);
\end{verbatim}
\end{data}

\begin{data}\label{data:Fermat-S4} The $S_4$-action on the 27 lines of the Fermat cubic surface, given by permuting coordinates on $\CP^3$, are generated by the following transposition and 4-cycle:
\begin{center}
    \begin{tabular}{l | l}
    \textbf{elt} & \textbf{permutation} \\
    \hline
    transp. & \texttt{(3,4)(5,11)(6,12)(7,9)(8,10)(13,15)(14,16)(17,21)(18,23)(19,22)(20,24)(26,27)} \\
    4-cycle & \texttt{(1,11,3,10)(2,12,4,9)(5,8,6,7)(13,23)(14,24,
    15,21)(16,22)(17,18,20,19)(25,27)} 
    \end{tabular}
\end{center}

\end{data}

\begin{data}\label{data:gens-K4-x-K4} 
The generators $\sigma_1,\sigma_2,\tau_1,\tau_2 \in W(E_6)$ from \autoref{StabIntersectNormalizer} are given by the following permutations:
\begin{center}
    \begin{tabular}{l | l}
    \textbf{elt} & \textbf{permutation} \\
    \hline
    $\sigma_1$ & \texttt{(1,3)(2,4)(5,6)(7,8)(9,12)(10,11)(14,15)(17,20)(18,19)(21,24)} \\
    $\sigma_2$ & \texttt{(1,4)(2,3)(5,8)(6,7)(9,10)(11,12)(13,16)(17,20)(21,24)(22,23)} \\
    $\tau_1$ & \texttt{(13,23)(14,19)(15,18)(16,22)(17,24)(20,21)} \\
    $\tau_2$ & \texttt{(1,4)(2,3)(9,11)(10,12)(13,16)(22,23)}
    \end{tabular}
\end{center}
\end{data}

\begin{data}\label{data:symmetric-monodromy} The (non-identity) elements in the Klein 4-group corresponding to symmetric monodromy are given by
\begin{center}
\begin{tabular}{l | l}
\textbf{elt} & \textbf{permutation} \\
\hline
$\tau_1$ & \texttt{(13,23)(14,19)(15,18)(16,22)(17,24)(20,21)}\\
$\sigma_1 \tau_2$ & \texttt{(1,3)(2,4)(5,6)(7,8)(9,12)(10,11)(13,23)(14,18)(15,19)(16,22)(17,21)(20,24)}\\
$\sigma_1 \tau_1 \tau_2$ & \texttt{(1,2)(3,4)(5,6)(7,8)(9,10)(11,12)(13,22)(14,18)(15,19)(16,23)(17,21)(20,24)}\\
\end{tabular}
\end{center}
\end{data}


\bibliography{symcubic}
\bibliographystyle{amsalpha}

\end{document}

%% file: commands.tex
\newcommand{\ZZ}{\mathbb{Z}}
\newcommand{\PP}{\mathbb{P}}
\newcommand{\FF}{\mathbb{F}}
\newcommand{\DD}{\mathbb{D}}

\newcommand{\CC}{\mathbb{C}}
\newcommand{\HH}{\mathbb{H}}

\newcommand{\CP}{\mathbb{CP}}

\DeclareMathOperator{\Aut}{Aut}

\DeclareMathOperator{\GL}{GL}
\DeclareMathOperator{\Gr}{Gr}

\DeclareMathOperator{\PGL}{PGL}

\DeclareMathOperator{\SL}{SL}

\DeclareMathOperator{\Stab}{Stab}

\DeclareMathOperator{\PO}{PO}

\providecommand{\PGL}{\text{PGL}}

\usepackage{xcolor}

\definecolor{purple}{RGB}{128, 0, 255}

\definecolor{darkgreen}{rgb}{0,0.30,0}

\definecolor{darkred}{rgb}{0.75,0,0}

\definecolor{darkblue}{rgb}{0,0,0.6}

\let\minus\smallsetminus

%% file: thms.tex
\usepackage{cleveref}

\def\makeautorefname#1#2{\expandafter\def\csname#1autorefname\endcsname{#2}}

\makeautorefname{eqn}{Equation}%
\makeautorefname{sec}{Section}%
\makeautorefname{subsec}{Subsection}%
\makeautorefname{footnote}{footnote}%
\makeautorefname{item}{item}%
\makeautorefname{figure}{Figure}%
\makeautorefname{table}{Table}%
\makeautorefname{wraptab}{wraptable}%
\makeautorefname{part}{Part}%
\makeautorefname{app}{Appendix}%
\makeautorefname{cla}{claim}%
\makeautorefname{assump}{assumption}%
\makeautorefname{cond}{conditions}%
\makeautorefname{conj}{conjecture}%
\makeautorefname{conv}{convention}%
\makeautorefname{cor}{corollary}%
\makeautorefname{cex}{counterexample}%
\makeautorefname{cexs}{counterexamples}%
\makeautorefname{data}{data}%
\makeautorefname{dig}{digression}%
\makeautorefname{disc}{discussion}%
\makeautorefname{def}{definition}%
\makeautorefname{ex}{example}%
\makeautorefname{exs}{examples}%
\makeautorefname{fac}{fact}%
\makeautorefname{intu}{intuition}%
\makeautorefname{lem}{lemma}%
\makeautorefname{meta}{metathm}%
\makeautorefname{nota}{notation}%
\makeautorefname{note}{note}%
\makeautorefname{warn}{warn}%
\makeautorefname{prop}{proposition}%
\makeautorefname{rmk}{remark}%
\makeautorefname{term}{terminology}%
\makeautorefname{thm}{theorem}%
\makeautorefname{upsh}{upshot}%
\theoremstyle{definition}
\newtheorem{theorem}{Theorem}[section]

\newtheorem{computation}[theorem]{Computation}

\newtheorem{convention}[theorem]{Convention}
\newtheorem{corollary}[theorem]{Corollary}

\newtheorem{data}[theorem]{Data}
\newtheorem{definition}[theorem]{Definition}

\newtheorem{example}[theorem]{Example}

\newtheorem{lemma}[theorem]{Lemma}

\newtheorem{proposition}[theorem]{Proposition}

\newtheorem{problem}[theorem]{Problem}
\newtheorem{remark}[theorem]{Remark}

\makeatletter
\let\c@corollary=\c@theorem
\let\c@proposition=\c@theorem
\let\c@lemma=\c@theorem
\let\c@assumption=\c@theorem
\let\c@conjecture=\c@theorem
\let\c@definition=\c@theorem
\let\c@example=\c@theorem
\let\c@remark=\c@theorem
\let\c@notation=\c@theorem
\let\c@equation\c@theorem
\makeatother

%% file: symcubic.bib
@article{Luna,
  title={Slices {\'e}tales},
  author={Luna, Domingo},
  journal={Sur les groupes alg{\'e}briques},
  volume={33},
  pages={81--105},
  year={1973}
}

@incollection{PopovVinberg,
  title={Invariant theory},
  author={Popov, Vladimir L and Vinberg, Ernest B},
  booktitle={Algebraic Geometry IV: Linear Algebraic Groups Invariant Theory},
  pages={123--278},
  year={1994},
  publisher={Springer}
}

@article{Sakamaki,
  title={Automorphism groups on normal singular cubic surfaces with no parameters},
  author={Sakamaki, Yoshiyuki},
  journal={Transactions of the American Mathematical Society},
  volume={362},
  number={5},
  pages={2641--2666},
  year={2010}
}

@article {Manivel,
    AUTHOR = {Manivel, L.},
     TITLE = {Configurations of lines and models of {L}ie algebras},
   JOURNAL = {J. Algebra},
  FJOURNAL = {Journal of Algebra},
    VOLUME = {304},
      YEAR = {2006},
    NUMBER = {1},
     PAGES = {457--486},
      ISSN = {0021-8693,1090-266X},
   MRCLASS = {17B25 (14J26 51E10 51E15 51E20)},
  MRNUMBER = {2256401},
MRREVIEWER = {Dmitry\ A.\ Timash\"ev},
       DOI = {10.1016/j.jalgebra.2006.04.029},
       URL = {https://doi.org/10.1016/j.jalgebra.2006.04.029},
}

@misc{EEG,
      title={Equivariant enumerative geometry}, 
      author={Thomas Brazelton},
      year={2024},
      eprint={2210.08622},
      archivePrefix={arXiv},
      primaryClass={math.AT},
      url={https://arxiv.org/abs/2210.08622}, 
}

@book{Margulis,
  title={{Discrete subgroups of semisimple Lie groups}},
  author={Margulis, Gregori A},
  volume={17},
  year={1991},
  publisher={Springer Science \& Business Media}
}

@article{FarbHandel,
  title={{Commensurations of $\mathrm{Out}(F_n) $}},
  author={Farb, Benson and Handel, Michael},
  journal={{Publications Math{\'e}matiques de l'IH{\'E}S}},
  volume={105},
  pages={1--48},
  year={2007}
}

@article{grothendieck1965elements,
  title={{{\'E}l{\'e}ments de g{\'e}om{\'e}trie alg{\'e}brique: IV. {\'E}tude locale des sch{\'e}mas et des morphismes de sch{\'e}mas, Seconde partie}},
  author={Grothendieck, Alexander},
  journal={Publications Math{\'e}matiques de l'IH{\'E}S},
  volume={24},
  pages={5--231},
  year={1965}
}

@article{hilbertVollen,
  title={{{\"U}ber die vollen Invariantensysteme}},
  author={Hilbert, David},
  journal={Math. Ann.},
  volume={42(3)},  
  pages={313–373},
  year={1893}
}

@article{MaclachlanHarvey,
  title={{On mapping-class groups and Teichm{\"u}ller spaces}},
  author={Maclachlan, Colin and Harvey, William J},
  journal={Proceedings of the London Mathematical Society},
  volume={3},
  number={4},
  pages={496--512},
  year={1975},
  publisher={Oxford University Press}
}

@article{BirmanHilden,
  title={{On isotopies of homeomorphisms of Riemann surfaces}},
  author={Birman, Joan S and Hilden, Hugh M},
  journal={Annals of Mathematics},
  volume={97},
  number={3},
  pages={424--439},
  year={1973},
  publisher={JSTOR}
}

@article{FarbWeinberger,
  title={Isometries, rigidity and universal covers},
  author={Farb, Benson and Weinberger, Shmuel},
  journal={Annals of Mathematics},
  pages={915--940},
  year={2008},
  publisher={JSTOR}
}

@inproceedings{RealACT,
  title={Hyperbolic geometry and moduli of real cubic surfaces},
  author={Allcock, Daniel and Carlson, James A and Toledo, Domingo},
  booktitle={Annales scientifiques de l'Ecole normale sup{\'e}rieure},
  volume={43},
  number={1},
  pages={69--115},
  year={2010}
}

@article{LandesmanLittSawin,
  title={{Big monodromy for higher Prym representations}},
  author={Landesman, Aaron and Litt, Daniel and Sawin, Will},
  journal={arXiv preprint arXiv:2401.13906},
  year={2024}
}

@article{ZhengOrbifold,
  title={Orbifold aspects of certain occult period maps},
  author={Zheng, Zhiwei},
  journal={Nagoya Mathematical Journal},
  volume={243},
  pages={137--156},
  year={2021},
  publisher={Cambridge University Press}
}

@article{YuZheng,
  title={Moduli spaces of symmetric cubic fourfolds and locally symmetric varieties},
  author={Yu, Chenglong and Zheng, Zhiwei},
  journal={Algebra \& Number Theory},
  volume={14},
  number={10},
  pages={2647--2683},
  year={2020},
  publisher={Mathematical Sciences Publishers}
}

@article{AllcockCarlsonToledo,
  title={The complex hyperbolic geometry of the moduli space of cubic surfaces},
  author={Allcock, Daniel and Carlson, James and Toledo, Domingo},
  journal={Journal of Algebraic Geometry},
  volume={11},
  number={4},
  pages={659--724},
  year={2002}
}

@book{Jordan,
  title={{Trait{\'e} des substitutions et des {\'e}quations alg{\'e}briques}},
  author={Jordan, Camille},
  year={1870},
  publisher={Gauthier Villars}
}

@article{Harris,
  title={Galois groups of enumerative problems},
  author={Harris, Joe},
  journal={Duke Math. J.},
  volume={46},
  number={1},
  pages={685--724},
  year={1979}
}

@book{HusemollerMilnor,
  title={Symmetric bilinear forms},
  author={Milnor, John Willard and Husemoller, Dale},
  volume={73},
  year={1973},
  publisher={Springer}
}

@misc{TheATLAS,
  title={{ATLAS of Finite Groups}},
  author={Conway, JH and Curtis, RT and Norton, SP and Parker, RA and Wilson, RA},
  year={1985},
  publisher={Clarendon Press, Oxford}
}

@article{LeykinSottile,
  title={{Galois groups of Schubert problems via homotopy computation}},
  author={Leykin, Anton and Sottile, Frank},
  journal={Mathematics of Computation},
  volume={78},
  number={267},
  pages={1749--1765},
  year={2009}
}

@article{SottileYahl,
  title={Galois groups in enumerative geometry and applications},
  author={Sottile, Frank and Yahl, Thomas},
  journal={arXiv preprint arXiv:2108.07905},
  year={2021}
}

@article{BeauvilleBourbaki,
  title={{Moduli of cubic surfaces and Hodge theory (after Allcock, Carlson, Toledo)}},
  author={Beauville, Arnaud},
  journal={G{\'e}om{\'e}tries {\`a} courbure n{\'e}gative ou nulle, groupes discrets et rigidit{\'e}s},
  volume={18},
  pages={445--466},
  year={2009},
  publisher={Soci{\'e}t{\'e} Math{\'e}matique de France Paris}
}

@inproceedings{BeauvilleMonodromy,
  title={Le groupe de monodromie des familles universelles d'hypersurfaces et d'intersections compl{\`e}tes},
  author={Beauville, Arnaud},
  booktitle={Complex Analysis and Algebraic Geometry: Proceedings of a Conference held in G{\"o}ttingen, June 25--July 2, 1985},
  pages={8--18},
  year={2006},
  organization={Springer}
}

@article{KudlaRapoport,
  title={On occult period maps},
  author={Kudla, Stephen and Rapoport, Michael},
  journal={Pacific Journal of Mathematics},
  volume={260},
  number={2},
  pages={565--582},
  year={2012},
  publisher={Mathematical Sciences Publishers}
}

@article{BruceWall,
  title={On the classification of cubic surfaces},
  author={Bruce, James W and Wall, Charles Terence Clegg},
  journal={Journal of the London Mathematical Society},
  volume={2},
  number={2},
  pages={245--256},
  year={1979},
  publisher={Oxford University Press}
}

@article{Allcock,
  title={New complex-and quaternion-hyperbolic reflection groups},
  author={Allcock, Daniel},
  journal={Duke Math. J.},
  volume={104},
  number={1},
  pages={303--333},
  year={2000}
}

@article{GriffithsRational,
  title={{On the periods of certain rational integrals: I, II}},
  author={Griffiths, Phillip A},
  journal={Annals of Mathematics},
  volume={90},
  number={2,3},
  pages={460--541},
  year={1969},
  publisher={JSTOR}
}

@misc{DuffLee,
      title={{Certified homotopy tracking using the Krawczyk method}}, 
      author={Timothy Duff and Kisun Lee},
      year={2024},
      eprint={2402.07053},
      archivePrefix={arXiv},
      primaryClass={math.NA},
      url={https://arxiv.org/abs/2402.07053}, 
}

@misc{Pandora,
    author = {Taylor Brysiewicz},
    title = {PANDO(RA): Parallel, Automated, Numerical, Discovery and Optimization (Research Aid)},
    howpublished ={https://github.com/tbrysiewicz/Pandora},
    year={2024}

}

@misc{Alberto,
    author = {Alberto Landi},
title = {Stacks, monodromy and symmetric cubic surfaces},
year = {2025},
publisher={Preprint}
}

@manual{GAP4,
    organization = "The GAP~Group",
    title        = "{GAP -- Groups, Algorithms, and Programming,
                    Version 4.14.0}",
    year         = 2024,
    url          = "\url{https://www.gap-system.org}",
    }

@article {three_skew,
    AUTHOR = {McKean, Stephen and Minahan, Daniel and Zhang, Tianyi},
     TITLE = {All lines on a smooth cubic surface in terms of three skew
              lines},
   JOURNAL = {New York J. Math.},
  FJOURNAL = {New York Journal of Mathematics},
    VOLUME = {27},
      YEAR = {2021},
     PAGES = {1305--1327},
      ISSN = {1076-9803},
   MRCLASS = {14N15 (14J26 14J70)},
  MRNUMBER = {4312735},
MRREVIEWER = {I.\ Dolgachev},
}

@misc{pichonpharabod2025galoisgroupssymmetriccubic,
      title={Galois Groups of Symmetric Cubic Surfaces}, 
      author={Eric Pichon-Pharabod and Simon Telen},
      year={2025},
      eprint={2509.06785},
      archivePrefix={arXiv},
      primaryClass={math.AG},
      url={https://arxiv.org/abs/2509.06785}, 
}
